\documentclass[twocolumn]{autart}

\usepackage{graphicx}          
\usepackage{natbib}        
\usepackage[utf8]{inputenc}
\usepackage[T1]{fontenc}
\usepackage{amsmath}
\usepackage{amsfonts}
\usepackage{amssymb}
\usepackage{algorithm}
\usepackage{algpseudocode}
\usepackage{caption}
\usepackage{subcaption}
\usepackage[usenames, dvipsnames]{color}
\usepackage{xcolor}
\usepackage{mathtools}
\usepackage{pgfplots}
\usepackage{wrapfig}

\newtheorem{corollary}{Corollary}
\newtheorem{assumption}{Assumption} 
\newtheorem{lemma}[thm]{Lemma}
\newtheorem{proposition}[thm]{Proposition}

\newtheorem{remark}{Remark}
\newtheorem{example}{Example}
\newtheorem{definition}{Definition} 
\AtBeginDocument{}

\begin{document}

\begin{frontmatter}

\title{Heuristic search for linear positive systems\thanksref{footnoteinfo}}

\thanks[footnoteinfo]{The authors are with the ELLIIT Strategic Research Area at Lund University, Lund, Sweden. This work is partially funded by Wallenberg AI, Autonomous Systems and Software Program (WASP) funded by the Knut and Alice Wallenberg Foundation, and the European Research Council (ERC) under the European Union’s Horizon 2020 Research and Innovation Programme under grant agreement No 834142 (ScalableControl).}

\author[First]{David Ohlin}\ead{david.ohlin@control.lth.se},
\author[First]{Anders Rantzer}\ead{anders.rantzer@control.lth.se},
\author[First]{Emma Tegling}\ead{emma.tegling@control.lth.se}

\address[First]{Department of Automatic Control, Lund University, P.O. Box 118, SE-221 00 Lund, Sweden}


\begin{abstract} 
This work considers infinite-horizon optimal control of positive linear systems applied to the case of network routing problems. We demonstrate the equivalence between Stochastic Shortest Path (SSP) problems and optimal control of a certain class of linear systems. This is used to construct a heuristic search framework for linear {positive} systems inspired by existing methods for SSP. {We propose a heuristics-based algorithm for {efficiently} finding local solutions to the analyzed class of optimal control problems with {a given initial} state and {positive} linear dynamics.} {By leveraging the bound on optimality in each state provided by the heuristics, we also derive a novel distributed algorithm for calculating local controllers within a specified performance bound, with a distributed condition for termination.} More fundamentally, the results allow for analysis of the conditions for explicit solutions to the Bellman equation utilized by heuristic search methods.
\end{abstract}

\end{frontmatter}

\section{Introduction}

The work presented here demonstrates the equivalence of two problem classes that are the focus of optimal control and graph search algorithms respectively. The first is a class of {optimal control problems involving} positive linear systems with linear cost and homogeneous constraints. The second is the class of Stochastic Shortest Path (SSP) problems and their generalizations, which are the focus of contemporary work on heuristic search algorithms, see \citep{bonet10heuristic} for a summary. In bridging the gap between these two formulations, we show that many results derived for either can be readily applied to the other. 

Specifically, we give a framework for a heuristics-based solution of optimal control problems with continuous state and action space that belong to the above described class of positive systems. In a similar spirit to heuristic search algorithms developed for SSP, these methods are computationally feasible for high-dimensional state and action spaces. {Additionally, we show that value iteration using both upper and lower heuristic bounds can leverage this measure of suboptimality to decentralize calculation with guarantees on global performance.} The analogy between the two classes also gives an analytic framework for further development of existing methods for SSP using results within the theory of optimal control. 

The method of dynamic programming, developed in \citep{bellman57dynamic}, shares a common origin with Dijkstra's algorithm, see \citep{dijkstra59note}. As emphasized in \citep{sniedovich06revisited}, the line of argument in Dijkstra's original work closely mirrors the optimality principle as first formulated by Bellman. 
Since their inception, these methods have both gone on to inspire a wide range of extensions, seeing great success in their respective fields of control theory and computer science. Although the lines of work within each discipline may seem disparate, they share the common feature of utilizing explicit solutions to the Bellman equation to achieve computational efficiency. It is this feature that enables the present work to give favorable results for positive systems and interpret these in terms of their counterpart in the SSP formulation.

A natural extension to Dijkstra's algorithm is proposed in \citep*{hart68heuristic}. The authors introduce a heuristic function that lower bounds the cost in any state and uses this to guide the search, prioritizing directions that give a low heuristic cost. Already in \citep{bertsekas91stochastic} the theoretical results of dynamic programming obtained in the deterministic case are extended to stochastic problems. Algorithms that implement similar ideas beyond the deterministic case through the use of heuristics are the subject of subsequent work in the field, see the seminal paper \citep*{barto95rtdp}. The close underlying connection between the study of heuristic functions and Lyapunov analysis is highlighted in \citep{perkins01lyapunov}. Comprehensive theoretical results on control of the broader class of constrained Markov Decision Processes (MDP), which includes problems like SSP, can be found in the monograph \citep{altman99cmdp}. Therein, arguments based in dynamic programming and Lyapunov theory are used to formulate linear programs for problems in this class under certain conditions.

In \citep*{mcmahan05brtdp} the authors incorporate the use of a heuristic upper bound and provide methods to find and improve such a bound. The extension is motivated by computational advantages, since the total cost is known to be upper bounded from the problem formulation of SSP. In this paper, the upper bound will additionally serve the purpose of guaranteeing the stability of solutions, which is not otherwise certain in the setting of optimal control. Experimental results demonstrating the robust performance of heuristic search methods for problems that deviate from the prescribed structure can be found in \citep{mcmahan05dijkstra}, together with an extended discussion of priority calculations for this family of algorithms. \citep{smith06focused} seek to improve existing algorithms by using upper and lower bound heuristics to inform the search. The class of problems for which heuristic methods can be applied is expanded in \citep*{kolobov11generalized} by the introduction of a broader class of SSP and a corresponding algorithm. The {work \citep{schmoll18routing} builds on these previous results, further improving} existing methods for certain classes of problems that can be rewritten as SSP.

Within the field of optimal control, recent works like \citep*{blanchini23exact} and \citep{rantzer22explicit} focus on explicit solutions to the Bellman equation for subclasses of positive linear systems. {The latter work formulates a linear program for the synthesis of the optimal controller. Previously, it has been shown in \citep{aitrami07synthesis} that stabilizing controller gains for linear positive systems can be found via linear programming.} In \citep*{ohlin23optimal}, the class of applicable systems is extended to include coupled homogeneous input constraints, allowing for the modeling of network routing problems while maintaining these computational advantages. Similar properties of solutions to the Bellman equation for a set of MDP with cost function based on the Kullback-Leibler divergence are noted in \citep{todorov06linearly}, with applications to optimal transport on networks. The use of upper and lower bounds on a general cost function for optimal control problems is explored in \citep{lincoln06relaxing} as a means to relax the problem formulation, allowing solutions to be suboptimal to an extent, in exchange for improved performance. 

This work shows an analogy between a subset of the class of positive linear systems described in \citep{ohlin23optimal} and SSP, with solutions for one class being provably optimal for the other. We develop linear upper and lower bounds for the class of positive systems mentioned above, inspired by heuristic search algorithms, that also yield explicit solutions to the Bellman equation. Using these bounds, we propose an algorithm inspired by existing heuristic search methods for SSP to find local solutions to the described class of optimal control problems.

The remainder of the paper is structured as follows; Section 2 defines the class of optimal control problems analyzed in subsequent sections. The optimal cost and corresponding policy are characterized through the application of Dynamic Programming. In Section 3, an equivalence between the analyzed problem and the class of SSP is shown. A constructive proof of the equivalence is given, with examples to illustrate the procedure. Linear heuristic bounds are introduced in Section 4. Section 5 proposes a heuristic search algorithm for positive systems. Section 6 contains numerical examples illustrating the performance of the proposed algorithm. Finally, conclusions are presented in Section 7.
\vspace*{-2mm}
\subsection{Notation} 
Inequalities are applied element-wise for matrices and vectors throughout. When applied to a set, the notation~${|\cdot|}$ is taken to mean the cardinality of that set. Further, the notation $\mathbb{R}^n_+$  is used to denote the closed positive orthant of dimension $n$. The operator~$\min\{{v},0\}$ extracts the minimum element of {a vector $v$}, yielding zero if ${v}$ has no negative elements. Let $A^{-\top}$ denote the inverse transpose of the matrix $A$. The expressions $\mathbf{1}_{p\times q}$ and $\mathbf{0}_{p\times q}$ signify a matrix of ones or zeros, respectively, of the indicated dimension, with subscript omitted when the size is clear from the context. If the dimension is zero, this is to be interpreted as the empty matrix. {The matrix $A$ is said to be row or column substochastic if it holds that ${A\mathbf{1} \le \mathbf{1}}$ or ${\mathbf{1}^\top A \le \mathbf{1}^\top}$ respectively.}
\vspace*{-2mm}
\section{Problem setup}
\vspace*{-2mm}
Consider the infinite-horizon optimal control problem
\vspace*{-2mm}
\begin{equation}
    \label{eq:optprob}
    \begin{aligned}
        \textnormal{Minimize} &\;\;\; \sum\limits_{t=0}^{\infty}\left[ s^\top x(t) + r^\top u(t) \right] \;\textnormal{over}\; \{u(t)\}^\infty_{t=0}\\
        \textnormal{subject to} &\;\;\; x(t+1) = Ax(t) + Bu(t)\\
        &\;\;\; u(t) \ge 0, \;\;\; x(0) = x_0\\
        & \;\;\; \begin{matrix} \mathbf{1}^\top u_1(t) & \le & E_1^\top x(t) \\ \vdots & & \vdots \\ \mathbf{1}^\top u_n(t) & \le & E_n^\top x(t) \end{matrix}
    \end{aligned}
\end{equation}
with ${x_0 \in \mathbb{R}^n_+}$. The input signal~${u\in\mathbb{R}^m}$ is partitioned into~$n$ subvectors~$u_i{\in\mathbb{R}^{m_i}}$ {(possibly ${m_i = 0}$, implying that the corresponding $B_i$ is empty)}, so that ${m = \sum_{i = 1}^{n} m_i}$. The {matrix} ${B=\begin{bmatrix}B_1 \cdots B_n\end{bmatrix}\in\mathbb{R}^{n\times m}}$ and ${A\in \mathbb{R}^{n\times n}}$ with ${B_i\in\mathbb{R}^{n\times m_i}}$ define the linear dynamics{, with each submatrix $B_i$ determining the effect of the inputs $u_i$ associated with the $i$th constraint}. Let the matrices~$B_i$ be further subdivided into individual columns as ${B_i = \begin{bmatrix} B_{i1} \cdots B_{im_i} \end{bmatrix}}$ with ${B_{ij}\in\mathbb{R}^n}$. The costs connected to the states and actions are ${s \in \mathbb{R}^{n}_+}$ with ${s > 0}$ and ${r \in \mathbb{R}^{m}_+}$ with ${r_i\in\mathbb{R}^{m_i}_+}$ following the partition of~$u$ and $r_{ij}$ denoting the $j$th element of $r_i$. The constraints on the input signal~$u(t)$ are given by the matrix ${E=\begin{bmatrix}E_1 \cdots E_n\end{bmatrix}^\top\in\mathbb{R}^{n\times n}_+}$. The special case of singular $E$ typically corresponds to degenerate formulations and will not be treated here. Define the set $\mathcal{V}$ to contain the indices $\{1,\ldots,n\}$. Let ${K = \left[K_1^\top \cdots K_n^\top\right]^\top}$ be a feedback matrix with ${K_i\in\mathbb{R}^{m_i\times n}_+}$ and define
\begin{equation*}
    \mathcal{K} = \{K: (\forall i\in\mathcal{V})\; \textbf{1}^\top K_i = E_i^\top  \;\textnormal{or}\; \textbf{1}^\top K_i = \mathbf{0}\}.
\end{equation*}
This set contains all feedback matrices that give full or zero actuation of the inputs $u_i$.

\begin{remark}
    \textnormal{The above formulation with strictly positive $s$ ensures that the system is {observable (only $x = 0$ yields zero immediate cost)} with the specified cost function, meaning that the optimal feedback law must also be stabilizing {if the total cost is finite}. A more intricate formulation of this criterion is given in \citep{li23exact}. In the present setting, however, any {observable} positive system on the form \eqref{eq:optprob} can be transformed to an equivalent problem satisfying $s > 0$ by reallocation of the costs $s$ and $r$.}
\end{remark}

{Let the function $J^*(x)$ give the optimal cost of an optimization problem initialized at $x$. For a general problem with immediate cost $g(x(t),u(t))$ and dynamics ${x(t+1)=f(x(t),u(t))}$, the celebrated Bellman equation}
\begin{equation}\label{eq:bellman}
    {J^*(x) = \min_{u} \left[g(x,u) + J^*(f(x,u))\right]}
\end{equation}

{gives a condition to validate the optimality of $J^*$. The main tool used to find optimal solutions in the present proofs is the method of \textit{value iteration}, which corresponds to fixed point iteration of \eqref{eq:bellman}. Of particular use is the property that value iteration converges monotonically to the fixed point of \eqref{eq:bellman} under light assumptions on the initial value. For further reading, see \citep{bellman57dynamic}.}

\subsection{Positive networks with linear cost}

In anticipation of Theorem \ref{th:main} we introduce the following assumption, which guarantees that positivity of the state is preserved under the constraints of \eqref{eq:optprob}.
\begin{assumption}
\label{as:ABE}
    The matrices $A,B$ and the set $\mathcal{K}$ satisfy 
    \begin{equation*}
        (A+BK)x \ge 0
    \end{equation*}
    for all $K\in\mathcal{K}$ and all states $x\in\mathbb{R}^n_+$.
\end{assumption}
With this assumption we restate the following result, which gives an explicit expression for the optimal cost function and an equivalent linear program.
\begin{thm}{\citep{ohlin23optimal}}\label{th:main}
    Under Assumption \ref{as:ABE}, the following three statements are equivalent:
    \begin{itemize}
        \item[]($i$) The problem (\ref{eq:optprob}) has a finite value for $x_0 \in \mathbb{R}^{n}_+$.        
        \item[]($ii$) There exists $p \in \mathbb{R}^{n}_+$ satisfying the equation
        \begin{equation}
        \label{eq:p}
            p = s + A^\top p + \sum_{i = 1}^n \min \{ r_i +B_i^\top p, 0 \} E_i.
        \end{equation}
        \item[]($iii$) The value of the linear program
        \begin{align*}
            \textnormal{Maximize} &\;\;\; \mathbf{1}^\top p \;\textnormal{over}\; p \in \mathbb{R}^{n}_+\\
        \textnormal{subject to} &\;\;\; p \; {\le} \; s + A^\top p + \sum_{i = 1}^n z_i E_i\\
        &\;\;\; z_i \le r_{ij} + B_{ij}^\top p \;\;\; \textnormal{for} \;\;\; j = 1,...,m_i\\
        &\;\;\; z_i \le 0
        \end{align*}
        is bounded with solution $p$ fulfilling \eqref{eq:p}.
    \end{itemize}
    Given the existence of a $p$ as in ($ii$), the minimal value of~\eqref{eq:optprob} is $p^\top x_0$. The optimal linear state feedback law is then given by $u_i(t) = K_ix(t)$ with
    \begin{equation}\label{eq:K}
        K_i := \begin{bmatrix} \mathbf{0}_{j-1\times n} \\ E_i^\top \\ \mathbf{0}_{m_i-j\times n} \end{bmatrix} \;\textnormal{for}\;i=1,...,n,
    \end{equation}
    where the vector $E_i^\top$ enters at the $j$th row with $j$ being the index of the minimal element of $r_i+B_i^\top p$, provided it is negative. If all elements are nonnegative then ${K_i = \mathbf{0}_{m_i\times n}}$.
\end{thm}


\begin{remark}
    \textnormal{An alternative formulation of \eqref{eq:p} is}
    \begin{equation}\label{eq:Kform}
        p = s + A^\top p + \min_{K\in\mathcal{K}}\left[ K^\top(r +B^\top p) \right]
    \end{equation}
    \textnormal{where the dependence on the constraint matrices $E_i$ is built into the definition of the set $\mathcal{K}$. The optimal cost is attained by choosing $K$ as in \eqref{eq:K}. Although less explicit, this expression is more compact {and} will be used when convenient in the treatment of the heuristic search algorithm presented in Sections 4 and 5 to illustrate the use of different policies for subsets of the state space $\mathcal{V}$.}
\end{remark}

\section{The equivalent problem}

The extension of previous results on positive linear systems to the class of systems covered by Theorem \ref{th:main} lets us make a formal connection to other existing frameworks for network routing problems. As the definitions used in the literature concerning SSP (e.g. \citep{bertsekas91stochastic}, \citep{mcmahan05brtdp} and \citep{kolobov11generalized}) differ in interpretation, the following definition strives to be as broad as possible while maintaining applicability of the results in Theorem \ref{th:main}.

\begin{definition}[Stochastic Shortest Path Problem]\label{def:SSP}\hfill\\
    A Stochastic Shortest Path (SSP) problem is defined by the following properties:
    \begin{itemize}
        \item[1.] a finite set of discrete states $\mathcal{V}'$,
        \item[2.] an initial state $v_0\in\mathcal{V}'$,
        \item[3.] a non-empty set of absorbing goal states $\mathcal{G}\subset \mathcal{V}'$,
        \item[4.] a finite action set $\mathcal{A}(v)$ for each state $v$,
        \item[5.] a function $\mathcal{T}(v,a)$ giving a vector of transition probabilities from $v$ to all other states for action~${a\in\mathcal{A}(v)}$,
        \item[6.] a cost function $\mathcal{C}(v,a)$ associating each action with a nonnegative cost that satisfies $\mathcal{C}(v,a) = 0$ for $v\in\mathcal{G}$ and $\mathcal{C}(v,a) > 0$ for $v\in\mathcal{V}'\setminus\mathcal{G}$
        \item[7.] there exists a policy $\pi(v_0)\in\mathcal{A}(v_0)$ for all initial states $v_0$ that reaches $\mathcal{G}$ with probability 1,
    \end{itemize}
    \vspace*{-2mm}
    with the objective
    \begin{equation*}
        \textnormal{Minimize}\;\;\; \mathbb{E}\left(\sum_{t = 0}^{\infty} \mathcal{C}(v_t,a_t)\right).
    \end{equation*}
\end{definition}
\vspace*{-2mm}
Note that the requirement $\mathcal{C}(v,a) > 0$ for $v\in\mathcal{V}'\setminus\mathcal{G}$ is equivalent to assuming {observability} in the setting \eqref{eq:optprob}. \citep{kolobov11generalized} lifts this restriction to allow zero-cost cycles in states outside $\mathcal{G}$ to be part of the solution. Such cases require special treatment and will not be considered in the present work.

\subsection{Reformulation}\label{sec:construction}

The problem \eqref{eq:optprob} {contains the class of SSP}. In order to show this, we first introduce a helpful lemma, as well as an assumption on the system dynamics. The following lemma simplifies the subsequent problem construction by allowing the choice $E = I$ {when the original $E$ is invertible}.

\begin{lemma}\label{lem:EI}
    A problem on the form \eqref{eq:optprob} with $E$ invertible can be transformed to an equivalent problem with $E = I$ that fulfills the requirement on positivity in Assumption~\ref{as:ABE}.
\end{lemma}
\vspace*{-4mm}
\begin{pf}
    Consider a similarity transform $\hat{x} = Tx$ of the original state vector in \eqref{eq:optprob}. This gives rise to transformed constraints ${\mathbf{1}^\top u_i \le E_i^\top T^{-1} \hat{x}}$. Enforcing ${ET^{-1} = I}$ yields ${T = E}$, which is a valid choice since~$E$ is invertible. The original state $x$ is confined to $\mathbb{R}^n_+$ under the dynamics of~\eqref{eq:optprob}. Further, the matrix $E$ is nonnegative by construction. Hence, it holds that ${\hat{x}\ge0}$ starting from ${\hat{x}_0 = Tx_0}$.\hfill $\blacksquare$
\end{pf}
\vspace*{-4mm}
Next, we introduce an assumption on substochasticity of the closed loop system in \eqref{eq:optprob}. As is shown below, both classes can be interpreted as SSP. 

\begin{assumption}\label{as:stable}
    The matrix ${(A+BK)}$ is column substochastic for all $K\in\mathcal{K}$.
\end{assumption}

Without this assumption, the problem can only be reformulated as an equivalent SSP {under} a relaxation of Definition \ref{def:SSP}, to allow for an infinite extension of the state space{, as detailed below in Section 3.3}. The resulting problem can, however, still be formulated as an optimization over finitely many actions, and the search methods introduced in subsequent sections, inspired by methods for SSP, are equally applicable regardless of this condition.

We are now ready to construct an SSP according to Definition \ref{def:SSP} from a problem on the form \eqref{eq:optprob} (with $E = I$) which fulfills Assumption \ref{as:stable}, and will show that the solution to this SSP is also an optimal control law for \eqref{eq:optprob}, with the same expected cost. Each of the objects necessary for Definition \ref{def:SSP} are defined as follows:

\begin{figure}[t]
    \centering
    \includegraphics[width=.9\linewidth]{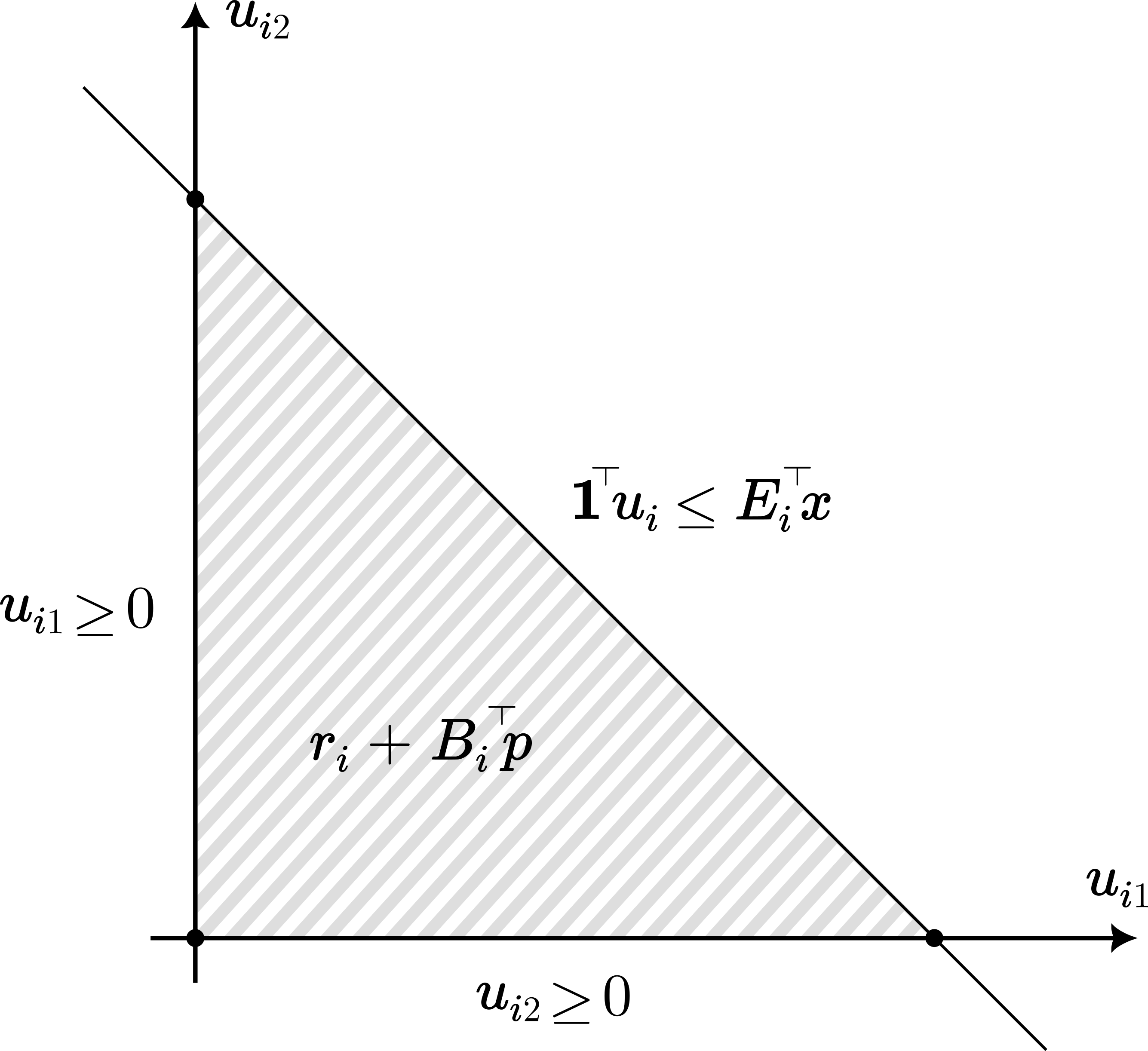}
    \caption{Allowed inputs under the constraints of \eqref{eq:optprob} for the case $m_i = 2$. The objective $(r_i+B_i^\top p)u_i$ is linear in the input $u_i$, so the optimum is attained on at least one of the vertices, for which an explicit expression is available as a linear function of the state $x$. A discrete action set $\mathcal{A}(i)$ can be constructed from the corresponding inputs.}
    \label{fig:constraint}
\end{figure}

\begin{itemize}
    \item[$\mathcal{V}'$:] Let each of the elements $x_v$, $v\in\mathcal{V}$ of the state vector in \eqref{eq:optprob} be identified with a corresponding state $v$ in a discrete set $\mathcal{V}_0$. Introduce a new, fictive state $v_g$ to act as an absorbing goal state. Define $\mathcal{V}' = \mathcal{V}_0 \cup \{v_g\}$.
    \item[$v_0$:] In order to have discrete initial states $v_0$, split~$x_0$ into components for each non-zero initial state. For each component $v$ the expected cost is that of the SSP starting from $v$, multiplied by the corresponding quantity in $x_0$.
    \item[$\mathcal{G}$:] The set $\mathcal{G}$ consists only of the fictive state $v_g$, which is absorbing by the definition of $\mathcal{A}(v_g)$ below.
    \item[$\mathcal{A}$:] A finite action space is created by noting that the optimum of the linear cost is attained (although not uniquely) for each subvector $u_v$ on the vertices of the constraint set given by ${\mathbf{1}\!^\top\!u_v(t)\le E_v^\top x}$ and ${u_v(t)\ge\mathbf{0}}$, according to Theorem~\ref{th:main} (see Figure~\ref{fig:constraint}). For each state these vertices generate a set of actions $\mathcal{A}(v)$.
    \item[$\mathcal{T}$:] The closed-loop system ${x(t+1) = (A+BK)x(t)}$ with $K\in\mathcal{K}$ is a transition from the state $x(t)$ encoded by the transition matrix $(A+BK)$. The resulting state $x(t+1)$ is interpreted as a probability distribution over the discrete states of $\mathcal{V}'$, and we define the transition function $\mathcal{T}(v,a)$ as the $v$th column of $(A+BK)$ where $K$ is given by the choice of action $a$. For all states $v$ and actions $a$ such that ${||\mathcal{T}(v,a)||_1<1}$, a transition to $v_g$ with probability ${1-||\mathcal{T}(v,a)||_1}$ is added. All that remains for the construction to be complete is to define the action set ${\mathcal{A}(v_g) = \{a_g\}}$, where $a_g$ gives unit transition probability to $v_g$.
    \item[$\mathcal{C}$:] Let $s_v$ denote the element of $s$ that multiplies $x_v$ and~$r_a$ denote the element of $r$ that multiplies the element of $u$ corresponding to the action $a\in\mathcal{A}(v)$. For actions~$a$ that correspond to ${u_v = \mathbf{0}}$ we let ${r_a = 0}$. The cost function is then naturally defined by associating each action with cost $\mathcal{C}(v,a) = s_v + r_a$. Finally, set ${\mathcal{C}(v_g,a_g) = 0}$.
\end{itemize}

\subsection{Equivalence}

Following this construction, we show equivalence in the sense of expected cost and optimal control law between the constructed SSP and the corresponding positive system treated in Theorem \ref{th:main}.

\begin{thm}\label{th:equivalence}
    Consider a problem on the form \eqref{eq:optprob} with finite value in which the system satisfies Assumptions \ref{as:ABE} and \ref{as:stable}. {Let $E$ be invertible. The optimal control problem} can be rewritten as an SSP according to Definition~\ref{def:SSP}, with equal expected cost. The optimal solution to this SSP is also an optimal control law for the problem~\eqref{eq:optprob}.
\end{thm} 
\vspace*{-4mm}
\begin{pf}
    By Lemma \ref{lem:EI}, it is sufficient to consider systems with $E = I$. Properties 1-6 of Definition \ref{def:SSP} are directly fulfilled by the construction in Section \ref{sec:construction}. All that remains is to show that {property} 7 also holds for the resulting SSP, that the expected cost is identical to that of the original problem and that the solution attains the optimal value of \eqref{eq:optprob}. 
    
    In order to show property 7 of Definition \ref{def:SSP} we apply the method used in Chapter 2.3 of \citep{bertsekas96neuro} to rewrite any discounted MDP as an SSP.  In the system ${x(t+1) = (A+BK)x(t)}$, the quantity~$||\mathcal{T}(v,a)||_1$ corresponds to the column sums of ${(A+BK)}$. The matrix ${(A+BK)}$ is Schur stable for some choices of ${K\in\mathcal{K}}$, since the problem \eqref{eq:optprob} has a finite value. In particular, this is at least true for the optimal feedback law. By the Perron-Frobenius theorem, we know that at least one column of a Schur stable~${(A+BK)}$ has sum strictly smaller than 1. Further, any policy $K$ that results in a stable system~${(A+BK)}$ will drive the state (excluding~$v_g$) to the origin. As~$v_g$ is absorbing, all probability density must eventually be concentrated in~$v_g$.
    
    The construction of $\mathcal{C}(v,a)$ in Section \ref{sec:construction} implies that the expected cost of some probability distribution~$x$ over~$\mathcal{V}'$ with a choice of actions $a\in\mathcal{A}(v)$ for each state is given by the sum $\sum_{v\in\mathcal{V}'} x_v \mathcal{C}(v,a)$, which is identical to the linear immediate cost of \eqref{eq:optprob}:
    \begin{align*}
        s^\top x + r^\top u &= s^\top x + r^\top Kx\\
        &= (s+\hat{r})^\top x\\
        &= \sum_{v\in\mathcal{V}'} x_v \mathcal{C}(v,a)
    \end{align*}
    where we use ${E = I}$, and ${\hat{r}\in\mathbb{R}^n}$ is the vector containing the element of $r$ (or zero if $u_i = 0$) corresponding to the action chosen for each state. Further, for some choice of actions ${a\in\mathcal{A}(v)}$ for all ${v\in\mathcal{V}'}$, the resulting transition matrix ${(A+BK)}$ attains the optimal cost of \eqref{eq:optprob}. This must then be the optimal cost also for the resulting SSP, since the cost is identical, with the corresponding actions giving a feedback law that solves~\eqref{eq:optprob}. \hfill $\blacksquare$
\end{pf}
\vspace*{-4mm}
\begin{example}
    \textnormal{Consider a problem on the form \eqref{eq:optprob} with $n = 3$, $m = 4$, $x_0 = \left[2 \; 0 \; 1\right]^\top$ and dynamics}
    \begin{align*}
        A = \begin{bmatrix}
            0.4 & 0 & 0 \\ 0 & 0.6 & 0 \\ 0.4 & 0.4 & 0.4
        \end{bmatrix} , \;\;\; 
        B = \begin{bmatrix}
            -0.4 & 0.3 & 0 & 0.2 \\ 0.4 & -0.6 & -0.5 & 0.2 \\ 0 & 0.3 & 0 & -0.4
        \end{bmatrix}
    \end{align*}
    \textnormal{with associated costs}
    \begin{equation*}
        s = \begin{bmatrix}
            1 & 1 & 1
        \end{bmatrix}^\top, \;\;\;
        r = \begin{bmatrix}
            1 & 1 & 1 & 1
        \end{bmatrix}^\top.
    \end{equation*}
    \textnormal{Constraints on the form \eqref{eq:optprob} with ${m_1 = 1}$, ${m_2 = 2}$, ${m_3 = 1}$ and ${E=I}$ ensure positivity of the dynamics, fulfilling Assumption \ref{as:ABE}. Following the construction in Section \ref{sec:construction}, we get a state space $\mathcal{V}'$ such that ${|\mathcal{V}'| = 4}$ where the fourth state acts as fictive goal state. Let $v$ be the state that corresponds to $x_2$. We then have ${|\mathcal{A}(v)| = 3}$ with actions $a_1$, $a_2$ and $a_3$ corresponding to ${u_{2} = \left[ \;0 \; 0 \;\right]^\top}$, ${u_{2} = \left[E_2^\top x \;\; 0\right]^\top}$ and ${u_{2} = \left[0 \;\; E_2^\top x\right]^\top}$ (the vertices in Figure~\ref{fig:constraint}). The cost function associated with these actions has value ${\mathcal{C}(v,a_1) = s_2 = 1}$, ${\mathcal{C}(v,a_2) = s_2 + r_{21} = 2}$ and ${\mathcal{C}(v,a_3) = s_2 + r_{22} = 2}$. The transition function is constructed from the second column of ${(A+BK_a)}$ for ${a\in\mathcal{A}(v)}$. Should the transition probabilities sum to less than one, the remainder is allocated to the transition from $v$ to the goal state. This gives }
    \begin{align*}
    \mathcal{T}(v,a_1) &= \left[ \; 0 \;\; 0.6 \;\; 0.4 \;\; 0 \; \right]^\top\\
    \mathcal{T}(v,a_2) &= \left[\; 0.3 \;\; 0 \;\; 0.7 \;\; 0 \;\right]^\top\\
    \mathcal{T}(v,a_3) &= \left[\; 0 \;\; 0.1 \;\; 0.4 \;\; 0.5 \;\right]^\top
    \end{align*}
    \textnormal{The optimal cost of the original problem starting from $x_0$ can be recovered by solving the constructed SSP and taking the first and third states as $v_0$, then adding together the resulting expected costs, weighted by the quantities 2 and 1 in $x_0$ respectively.}
\end{example}

To complete the equivalence, the following theorem states that the converse also holds. Due to the constructive nature of Theorem \ref{th:equivalence}, the reverse direction follows in a straightforward manner from the same reasoning.

\begin{thm}\label{th:converse}
    Any SSP can be rewritten as an optimal control problem on the form \eqref{eq:optprob}, with finite value equal to the expected cost of the SSP, and dynamics satisfying Assumptions \ref{as:ABE} and \ref{as:stable}. The solution to the SSP also attains the optimum of the problem \eqref{eq:optprob}.
\end{thm}
\vspace*{-4mm}
\begin{pf}
    Taking $x_0$ as the vector where ${x_{v_0} = 1}$ and zero in all other elements gives the initial state. Extract ${\mathcal{V} = \mathcal{V}'\setminus{G}}$ as the set of states in \eqref{eq:optprob}. The converse direction then follows directly from reversing the construction in Section \ref{sec:construction} and applying the proof of Theorem \ref{th:equivalence}. The optimal controller for \eqref{eq:optprob} may not be unique. However, the controller given by \eqref{eq:K} is always a solution to the SSP, attaining the same expected cost. \hfill $\blacksquare$
\end{pf}
\vspace*{-4mm}
To the best of our knowledge, the equivalence resulting from Theorems~\ref{th:equivalence} and~\ref{th:converse}, although intuitive, has not previously been observed in the literature. Despite this, there is much to be gained from making this connection explicit. Firstly, the success of methods developed to efficiently solve problems within the class of SSP can be utilized for optimal control of large-scale systems with continuous state and action spaces. The full class of problems on the form \eqref{eq:optprob}, i.e. without Assumption~\ref{as:stable}, is broader than the class of SSP as given in Definition~\ref{def:SSP}. However, as is exemplified by the heuristic search algorithm presented in Section 4, tools and algorithms developed for SSP can be applied to all problems on the form~\eqref{eq:optprob}, as the key advantages of the solution to the Bellman equation still hold. Secondly, the analysis provided in the optimal control framework provides a clear view of the problem class for which an explicit solution to the Bellman equation is obtained. This novel framework can be used to find previously unexplored classes of systems that share the computational advantages of SSP. One such example is the class of systems on the form of~\eqref{eq:optprob} with the summation constraints on the input vectors $u_i$ replaced by arbitrary $p$-norms. This setting is explored in~\citep{li23exact}, although for a more restricted problem class with decoupled constraints. 

\subsection{Extension}\label{sec:extension}

We have shown above that there exists an SSP formulation in accordance with Definition \ref{def:SSP} that is equivalent to \eqref{eq:optprob} when $(A+BK)$ is column substochastic for all ${K\in\mathcal{K}}$ (Assumption~2). By relaxing the first condition of Definition~\ref{def:SSP}, to allow~$\mathcal{V}'$ to contain an infinite number of states, it is possible to construct an equivalent SSP with unit sum transition probabilities, even when column sums of ${(A+BK)}$ exceed unity. 

\begin{figure}[t]
    \centering
    \includegraphics[width=\linewidth]{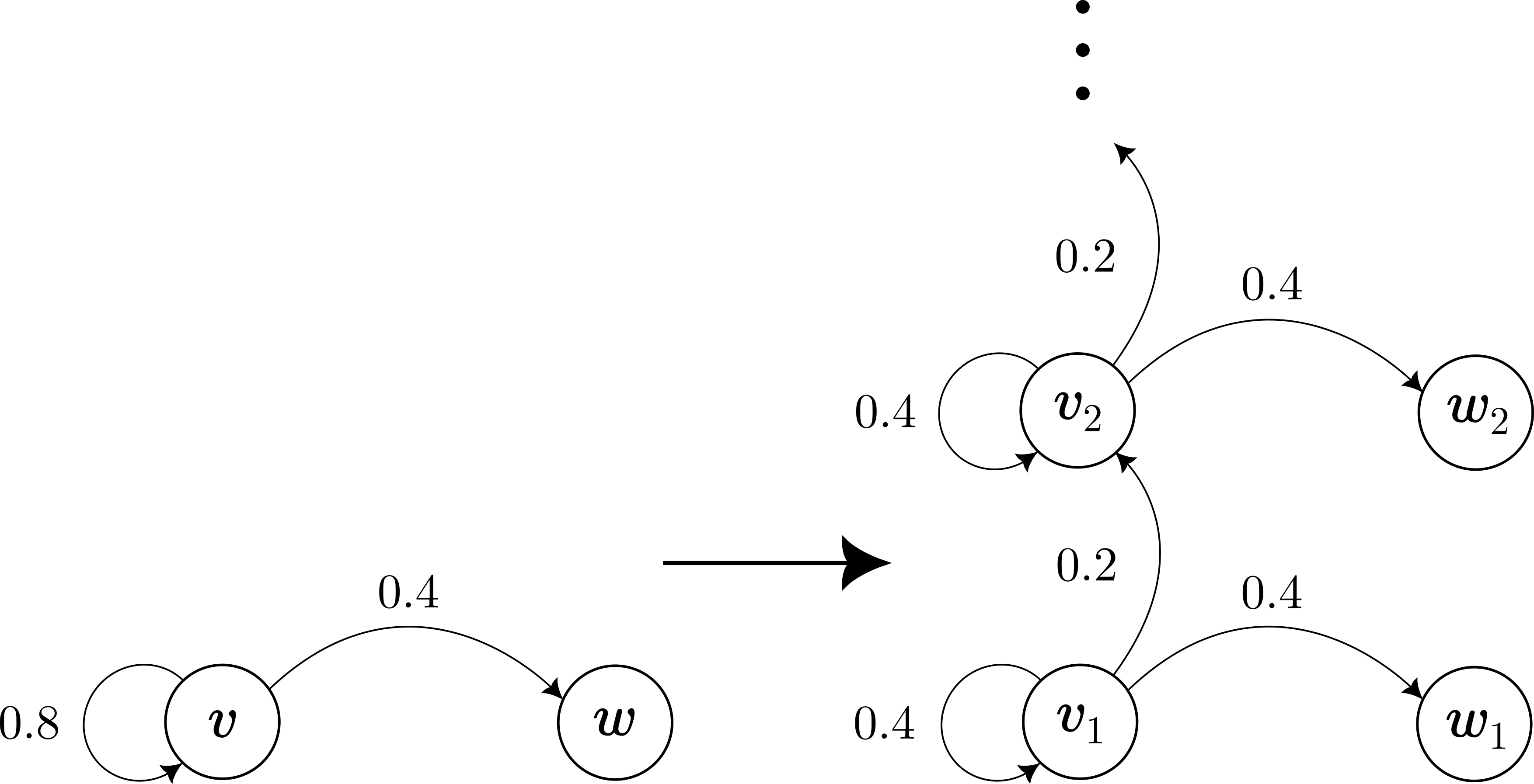}
    \caption{Expansion of the state space as described in Section~3.3 to accommodate the unit sum transitions of Definition~\ref{def:SSP}, with multiple states $v_k$ corresponding to the original $v$. A finite SSP cannot represent the unbounded growth in the state that is possible for, e.g., an unstable system in the formulation \eqref{eq:optprob}. {We must instead resort to an infinite space of artificial states with nearly identical properties. As is detailed in \citep{wrobel84skeleton}, there exist tractable methods for simplification of systems with this structure.}}
    \label{fig:expansion}
\end{figure}

Take an SSP generated by the method described in Section 3.1 that has ${||\mathcal{T}(v,a)||_1>1}$ for some~$v$, $a$. Let the transition associated with action $a\in\mathcal{A}(v)$ be defined by the probabilities $\{t_w(v,a)\}_{w\in\mathcal{V}'}$ so that action~$a$ gives the probability $t_w(v,a)$ of transitioning from state~$v$ to state~$w$ and $||\mathcal{T}(v,a)||_1 = \sum_{w\in\mathcal{V}'}t_w(v,a)$. For each state $v\in\mathcal{V}'$ generate a set ${\mathcal{V^+}(v) = \{v_k\}_{k=1}^{\infty}}$ where ${v_1 = v}$. Redefine the actions $a\in\mathcal{A}(v)$ for the original states $v\in\mathcal{V}'$ so that they fulfill the properties
\begin{equation}\label{eq:kt}
    \sum_{w_k\in\mathcal{V}^+(w)} kt_{w_k}(v,a) = t_w(v,a)
\end{equation}
and
\begin{equation}\label{eq:sumcond}
    \sum_{w\in\mathcal{V}'}\sum_{w_k\in\mathcal{V}^+(w)} t_{w_k}(v,a) = 1.
\end{equation}
The condition \eqref{eq:sumcond} is equivalent to ${||\mathcal{T}(v,a)||_1 = 1}$. For the states ${v_k\in\mathcal{V}^+(v)}$, ${k\ge2}$, define $\mathcal{A}(v_k)$ by taking each action ${a\in\mathcal{A}(v)}$ and creating a corresponding action ${a'\in\mathcal{A}(v_k)}$ such that
\begin{equation}\label{def:t}
    t_{w_{j}}(v_k,a') = \begin{cases}
        t_{w_{\ell}}(v_1,a) \;\;\; j = k\ell, \; \ell = 1,2,3...\\
        0 \;\;\;\;\;\;\;\;\;\;\;\;\;\;\; \textnormal{otherwise}
    \end{cases}
\end{equation}
and ${\mathcal{C}(v_k,a') = k\mathcal{C}(v_1,a)}$. The reallocation of transition probabilities in accordance with \eqref{eq:kt} and \eqref{eq:sumcond} is illustrated in Figure \ref{fig:expansion} for a single state. The result is formalized in the following proposition.

\begin{proposition}\label{prop:infty}
    If the set $\mathcal{V}'$ in Definition \ref{def:SSP} is not required to be finite, then Theorem \ref{th:equivalence} applies also to problems on the form~\eqref{eq:optprob} that violate Assumption \ref{as:stable}.
\end{proposition}
\vspace*{-4mm}
\begin{pf}
    The Bellman equation for the states of the new SSP becomes
    \begin{align*}
        J^*\!(v_k) &= \!\!\!\min_{a'\in\mathcal{A}(v_k)}\!\!\!\left[\!\mathcal{C}(v_k,a') \!+\!\!\!\sum_{w\in\mathcal{V}'}\!\!\!\!\!\!\!\!\sum_{\;\;\;\;\;\;w_j\in\mathcal{V}^+\!(w)}\!\!\!\!\!\!\!\!\!\!t_{w_j}\!(v_k,a')J^*\!(w_j)\!\right].
    \end{align*}
    We make the ansatz $J^*(v_k) = kJ^*(v_1)$. Applying \eqref{def:t} results in the equation
    \begin{align*}
        kJ^*\!(v_1) &= \!\!\!\!\min_{a\in\mathcal{A}(v_1)}\!\!\!\left[\!k\mathcal{C}(v_1,a) \!+\!\!\!\sum_{w\in\mathcal{V}'}\!\!\!\!\!\!\!\!\sum_{\;\;\;\;\;\;w_{\ell}\in\mathcal{V}^+\!(w)}\!\!\!\!\!\!\!\!\!\!t_{w_{\ell}}\!(v_1,a)J^*\!(w_{k\ell})\!\right]
    \end{align*}
    where, since $J^*(w_{k\ell}) = kJ^*(w_\ell)$, division by $k$ yields the Bellman equation for $v_1$, verifying the ansatz. To show that this is equivalent to the original problem, take the Bellman equation for the state $v_1$ in the new formulation and use the same ansatz, giving
    \begin{align*}
        J^*\!(v_1) &= \!\!\!\!\!\min_{a'\in\mathcal{A}(v_1)}\!\!\!\left[\mathcal{C}(v_1,a') \!+\!\!\!\sum_{w\in\mathcal{V}'}\!\!\!\!\!\!\!\!\sum_{\;\;\;\;\;\;w_k\in\mathcal{V}^+\!(w)}\!\!\!\!\!\!\!\!\!\!t_{w_k}\!(v_1,a')J^*\!(w_k)\right]\\
        &= \!\!\!\min_{a\in\mathcal{A}(v)}\!\left[\mathcal{C}(v,a) \!+\!\!\!\sum_{w\in\mathcal{V}'}\!\!\!\!\!\!\!\!\sum_{\;\;\;\;\;\;w_k\in\mathcal{V}^+\!(w)}\!\!\!\!\!\!\!\!\!\!t_{w_k}\!(v,a)kJ^*\!(w)\right]\\
        &= \!\!\!\min_{a\in\mathcal{A}(v)}\!\left[\mathcal{C}(v,a) \!+\!\!\!\sum_{w\in\mathcal{V}'} \!t_w(v,a)J^*\!(w)\right]
    \end{align*}
    where we have used \eqref{eq:sumcond}, and recall that ${v_1 = v}$. The final expression gives the optimal cost $J^*(v)$ of the initial problem, completing the proof.\hfill $\blacksquare$
\end{pf}
\vspace*{-4mm}

\begin{example}
    \textnormal{Take the system used in Example 1, but substitute the dynamics}
    \begin{equation*}
        A = \begin{bmatrix}
            0.4 & 0 & 0 \\ 0 & 0.8 & 0 \\ 0.4 & 0.4 & 0.4
        \end{bmatrix}
    \end{equation*}
\end{example}
\vspace*{-4mm}
\textnormal{Following the procedure in Section 3.1 results in the transition ${\mathcal{T}(v,a_1) = \left[ \; 0 \;\; 0.8 \;\; 0.4 \;\; 0 \; \right]^\top}$ for the first action considered in Example 1. By extending the state space as described in Section 3.3, we instead get transitions defined by $t_{v_1}(v_1,a_1) = 0.2$, $t_{v_2}(v_1,a_1) = 0.4$ and $t_{w_1}(v_1,a_1) = 0.4$ where $v_k$ and $w_k$ are the states corresponding to the second and third states after expansion. The situation is illustrated in Figure \ref{fig:expansion}.}

\textnormal{The infinite size of the state space as {required} in Proposition \ref{prop:infty} may seem discouraging at first glance{, but is necessary to translate the concept of instability to the setting of SSP}. However, the minimization over the action space for the extended states $v_k$ corresponding to some $v\in\mathcal{V}'$ in the original formulation is identical apart from scaling by a factor for all $k$. This means that the solution to the SSP can still be found by optimization over a finite set of actions, which is guaranteed to contain an optimal solution. Previous results in \citep{wrobel84skeleton} show the possibility of reducing infinite state MDP to a finite "skeleton" under similar conditions. For more discussion on SSP and heuristic search for infinite state spaces, see \citep{perkins01lyapunov} and \citep{bertsekas18proper}.}

\section{Heuristic bounds}

Instances of optimal control problems like (\ref{eq:optprob}) may have intractably large state and action spaces, prompting the search for methods that do not require calculations using full dynamics and state. In order to give concrete guarantees on performance, heuristic functions are introduced to bound the optimal cost in any given state. Within the area of control, this has taken the form of relaxed Dynamic Programming, as in \citep{lincoln06relaxing}. On the other hand, heuristic search methods that are a direct extension of Dijkstra's algorithm have been developed for deterministic shortest path problems, an early example being the $A^*$ algorithm presented in \citep{hart68heuristic}. These utilize physical intuition and assumptions about the structure and cost of the analyzed network to construct similar bounding heuristics that guide the search for an optimal solution. For the dynamical systems explored here, the prior knowledge of the system dynamics required to establish such bounds can be represented by some initial stabilizing policy{, or a policy stabilizing only a subset of the states.} We first define two functions to bound the optimal cost function given by the solution of \eqref{eq:p}.
\begin{definition}\label{def:h}
    The upper and lower bound heuristic functions are defined respectively as
    \begin{align*}
        \overline{H}(x) &= \overline{h}^\top x\\
        \underline{H}(x) &= \underline{h}^\top x.
    \end{align*}
\end{definition}




\subsection{Admissibility}

We next present an initialization of the heuristics, based on {a} stabilizing policy in the case of the upper bound. The following propositions shows {that \textit{admissibility} of} the heuristic bounds for this initial choice, and proves that these hold under value iteration.

\begin{proposition}\label{prop:h}
    Consider a system on the form \eqref{eq:optprob} that fulfills Theorem \ref{th:main}. Given {a stabilizing controller ${\hat{K}\in\mathcal{K}}$} and heuristic functions as in Definition \ref{def:h} with
    \begin{align*}
        \overline{h} &= (I - (A + B\hat{K}))^{-\top}(s+\hat{K}^\top r)\\
        \underline{h} &= s
    \end{align*}
    the following relation holds
    \begin{equation}\label{eq:bound}
        \underline{H}(x) \le J^*(x) \le \overline{H}(x).
    \end{equation}
    Further, value iteration 
    for either heuristic yields an explicit update on the form \eqref{eq:p} and preserves the relation~\eqref{eq:bound}.
\end{proposition}
\vspace*{-4mm}
\begin{pf}
    The immediate cost in any state $x(t)$ is given by~$g(x(t),u(t)) = s^\top x + r^\top u$. As a consequence, it trivially holds that $\underline{h}^\top x = s^\top x \le g(x(t),u(t)) \le J^*(x)$. For the second heuristic, we first inspect the expression given in the proposition. It is clear that, since the feedback law~$u = \hat{K}x$ stabilizes the system in \eqref{eq:optprob}, the eigenvalues $\lambda$ of the closed loop system fulfill $|\lambda(A+B\hat{K})| < 1$. This in turn implies that the inverse $(I-(A+B\hat{K}))^{-1}$ is well defined. Using the upper bound heuristic $\overline{H}(x)$, the associated cost of a state~$x$ is 
    \begin{align}
        \overline{h}^\top x &= (s+\hat{K}^\top r)^\top(I - (A + B\hat{K}))^{-1} x\nonumber\\
        \iff \overline{h}^\top x &= s^\top x + r^\top \hat{K} x + \overline{h}^\top(A+B\hat{K})x. \label{eq:kcost}
    \end{align}
    Next, consider the optimal cost ${J^*(x) = p^\top x}$, which is given by the equation
    \begin{align}
        p^\top x &=\!\min_{u\in U(x)}\!\left[g(x,u) + J^*(Ax + Bu)\right]\nonumber\\
        \iff\!p^\top\!x &= \min_{K\in\mathcal{K}} \!\left[s^\top x + r^\top Kx + p^\top(A+BK)x\right]\!. \label{eq:pcost}
    \end{align}
    Since we have $p,x\ge0$, ${p \le \overline{h}}$ follows from comparison of the expressions \eqref{eq:kcost} and \eqref{eq:pcost}. Equality is attained by the choice ${\hat{K}\in\mathcal{K}}$ in \eqref{eq:pcost}, directly giving the upper bound on~$p$. 
    
    Iteration of the Bellman equation yields an update on the form \eqref{eq:p} as a direct consequence of Theorem \ref{th:main}, since the heuristic function is linear and $\overline{h},\underline{h} \ge 0$. To show that the relation holds after iteration of \eqref{eq:p} it is sufficient to note that value iteration of $p \le \overline{h}{(k)}$ is a monotone operation with $p$ as a fixed point, together implying that $p\le\overline{h}{(k+1)}$ also holds. An identical argument can be made for the lower bound $\underline{h}$. \hfill $\blacksquare$
\end{pf}

\subsection{Consistency}

We have shown that heuristic bounds on the form proposed above indeed bound the optimal cost and can be improved through value iteration. However, it is not in general true for an arbitrary heuristic function that the bound in each state improves with every {fixed-point} iteration of \eqref{eq:p}. In the context of SSP, this property holds only for \textit{consistent} heuristics. The condition required for a lower bound heuristic to be consistent is, as has been noted elsewhere in the literature, equivalent to the Bellman inequality
\begin{equation*}
    \underline{H}(x) \le \min_{u\in U(x)}\left[g(x,u) + \underline{H}(Ax + Bu)\right]
\end{equation*}
\vspace*{-2mm}
which, for the class of problems \eqref{eq:optprob}, is equivalent to
\begin{equation}\label{eq:consistent}
    \underline{h} \le s + A^\top \underline{h} + \sum_{i = 1}^n \min \{ r_i +B_i^\top \underline{h}, 0 \} E_i.
\end{equation}
It is worth noting that this is the exact form of the constraints of the linear program ($iii$) in Theorem~\ref{th:main}, giving that program the natural interpretation of a search over consistent lower bound heuristics. For the upper bound heuristic $\overline{h}$, consistency instead corresponds to~\eqref{eq:consistent} with the inequality reversed. This property holds for the heuristics given below in Proposition \ref{prop:h} and is used to prove monotone improvement under value iteration.


{We can establish a finite upper heuristic bound in the case of an optimal control problem on the form \eqref{eq:optprob} if a stabilizing policy is known.} This allows for an upper bound of the cost in any given state, since following this (likely suboptimal) {policy} is guaranteed to stabilize the system, bringing the immediate cost $g(x(t),u(t))$ to zero asymptotically. This is necessary to guarantee finite cost {for all possible $x_0\in\mathbb{R}^n_+$} in the case of an unstable but stabilizable open-loop system. In the case of value iteration for the whole system this would not be a concern; the selected policy generating the minimal cost in each step must be stabilizing after a single iteration.

\begin{proposition}\label{prop:consistency}
    {Heuristic bound functions ${\underline{H}(x) = s}$ and $\overline{H}(x)$ corresponding to the cost under any fixed policy ${\hat{K}\in\mathcal{K}}$ are consistent, yielding} monotonically improving {value iteration} sequences
    \begin{equation*}
        \overline{h}{(0)} \ge \overline{h}{(1)} \ge \cdots \ge \overline{h}{(k)} \;\;\;\textnormal{and}\;\;\; \underline{h}{(0)} \le \underline{h}{(1)} \le \cdots \le \underline{h}{(k)}
    \end{equation*}
    where each inequality is strict in at least one element unless the fixed point is reached.
\end{proposition}
\vspace*{-4mm}
\begin{pf}
    Taking the definition of the initial upper bound heuristic in Proposition \ref{prop:h} as $\overline{h}{(0)}$ we get
    \begin{align*}
        \overline{h}{(0)} &= (I - (A + B\hat{K}))^{-\top}(s+\hat{K}^\top r)\\
         &= s + \hat{K}^\top r + (A+B\hat{K})^\top\overline{h}{(0)}\\
         &\ge s + A^\top\overline{h}{(0)} + \min_{K\in\mathcal{K}}\left[K^\top (r + B^\top\overline{h}{(0)})\right]\\
         &= \overline{h}{(1)}.
    \end{align*}
    Since the update \eqref{eq:p} is an instance of value iteration, which is a monotone operation, $\overline{h}{(0)}\ge\overline{h}{(1)}$ implies that $\overline{h}{(k)}\ge\overline{h}{(k+1)}$ for all $k$. For the lower bound we have $s > 0$. Further, taking $\underline{h}{(0)} = 0$ and applying \eqref{eq:p} yields
    \begin{align*}
        \underline{h}{(1)} &= s + A^\top\underline{h}{(0)} + \min_{K\in\mathcal{K}}\left[K^\top(r+B^\top\underline{h}{(0)})\right]\\
        &= s
    \end{align*}
    which is the value of the initial lower bound given in Proposition~\ref{prop:h}. Hence, $\underline{h}{(0)} \le \underline{h}{(1)}$ and the same argument as above gives monotone improvement for all $k$. Finally, strict improvement in at least one element follows directly from the fact that any cost function that is fixed under value iteration is also a solution to \eqref{eq:p}. \hfill $\blacksquare$
\end{pf}
\vspace*{-4mm}
Going forward, we assume that the value of the heuristic bounds are given a priori, as is customary in the literature on similar methods for SSP, {but not necessarily finite in all elements}. {When one or more states are unstable under an evaluated policy, we set the corresponding elements of the upper bound heuristic to $\overline{h}_i = \infty$, to represent the lack of available bounds.} This assumption may be taken to represent some prior theoretical knowledge, or physical intuition gained from practical experience of the system. However, it should be noted that, since the goal of methods using these heuristics is to avoid calculations scaling with the (potentially large) size of the system, it is not a trivial proposition to evaluate the given expression for~$\overline{h}$ {according to Proposition~\ref{prop:h}, even if some global} stabilizing controller $\hat{K}$ is known. {An analogous problem} in the setting of SSP is explored in \citep{mcmahan05brtdp}, wherein the authors present a framework for finding {and improving a globally} consistent heuristic bound (giving monotone improvement) that scales favorably with the size of the problem. 

\subsection{{Superconsistency}}

{Previous heuristic search methods for SSP give guarantees on monotone convergence given a consistent initial heuristic. However, this does not suffice for the upper bound heuristic in the full class of problems~\eqref{eq:optprob} unless a stabilizing policy is known. We therefore introduce the notion of \textit{superconsistency}. Let $E_{ii}$ denote the $i$th diagonal element of $E$. An upper bound heuristic $\overline{h}$ is superconsistent if it satisfies the following inequality:}

\begin{equation}\label{eq:super}
    {\overline{h} \ge s + A^\top \overline{h} + \begin{bmatrix}
        \min \{ r_1 +B_1^\top \overline{h}, 0 \} E_{11} \\ \vdots \\ \min \{ r_n +B_n^\top \overline{h}, 0 \} E_{nn}
    \end{bmatrix}}
\end{equation}

{The rightmost term is the vector corresponding to the sum $\sum_{i = 1}^n \min \{ r_i +B_i^\top p, 0 \} E_i$ in~\eqref{eq:p} with all off-diagonal elements of $E$ set to zero. As a consequence, it coincides with the normal notion of consistency for $E=I$, but is a stronger requirement in the general case. This property will be essential in proving monotone convergence for a general problem on the form~\eqref{eq:optprob} when an initial stabilizing policy is not known.}

\subsection{Rate of convergence}

The following convergence proof shows that the quality of the initial policy determines the rate of convergence for value iteration. This means that it is computationally desirable to start from accurate heuristics.
\begin{proposition}\label{prop:VI}
     Take some $\delta < 1$ such that $\underline{h}\ge\delta\overline{h}$. Further, let $\beta>1$ be some number that fulfills $A^\top\overline{h}\le\beta s$ and $B^\top\overline{h}\le\beta r$. After $k$ iterations of \eqref{eq:p} it holds that  
    \begin{align}
        \underline{h}{(k)}\ge\left(1 - \frac{1-\delta}{(1-\beta^{-1})^{k}}\right)p.
    \end{align}
    for the resulting heuristic $\underline{h}{(k)}$.
\end{proposition}
\vspace*{-4mm}
\begin{pf}
    The problem \eqref{eq:optprob} has a finite value for any system that fulfills Theorem \ref{th:main}, implying that $\overline{h}<\infty$. This guarantees the existence of $\beta$ and $\delta$ according to the given conditions. Applying the operator \eqref{eq:p} to the initial heuristic $\underline{h}{(0)}$ gives
    \begin{align*}
        \underline{h}{(1)} &= s + A^\top \underline{h}{(0)} + \sum_{i=1}^{M} \min \{ r_i +B_i^\top\underline{h}{(0)}, 0 \} E_i\\
        &\ge s + \delta A^\top p + \frac{1-\delta}{1+\beta}(\underbrace{A^\top p-\beta s}_{\le0})\\
        &+ \sum_{i=1}^{M} \min \{r_i + \delta B_i^\top p + \frac{1-\delta}{1+\beta}(\underbrace{B_i^\top p-\beta r_i}_{\le0}), 0 \} E_i\\
        &= \left(\!1\!-\!\frac{1-\delta}{(1\!-\!\beta^{-1})}\!\right)\!\underbrace{\left(\!s + A^\top\!p +\!\sum_{i=1}^{M}\!\min \{r_i\!+\!B_i^\top p, 0 \} E_i\!\right)}_{=\;p}
    \end{align*}
    The same calculation can be repeated to yield successive bounds up to iteration $k$. An identical bound for convergence of $\overline{h}$ follows from the same reasoning. \hfill $\blacksquare$
\end{pf}
\vspace*{-4mm}
The bound proven above depends on the global accuracy of the heuristics $\overline{h}$ and $\underline{h}$, indicating that local convergence is slow for an initial policy that performs poorly in some (possibly distant) state. In the following section we present a heuristics-based search algorithm that iteratively calculates local solutions without requirements on global optimality. As a consequence, the above statement can be used to show a faster rate of local convergence when accurate heuristics are available locally.

\section{Heuristic search for positive systems}

Next, we propose an algorithm that, starting from some initial state $x_0$, locally improves the heuristic bounds in Definition~\ref{def:h} by successively expanding a search area $S$. In keeping with the spirit of heuristic search algorithms for SSP, the aim of this method is to find a solution that is sufficiently close to optimality while avoiding an exhaustive search of the full state and action spaces. {The heuristic bounds serve both to measure the distance to optimality and to guide the search, giving faster execution the more detailed the initial information in $\overline{h}$ and $\underline{h}$.}

\begin{figure}[t]
    \centering
    \includegraphics[width=\linewidth]{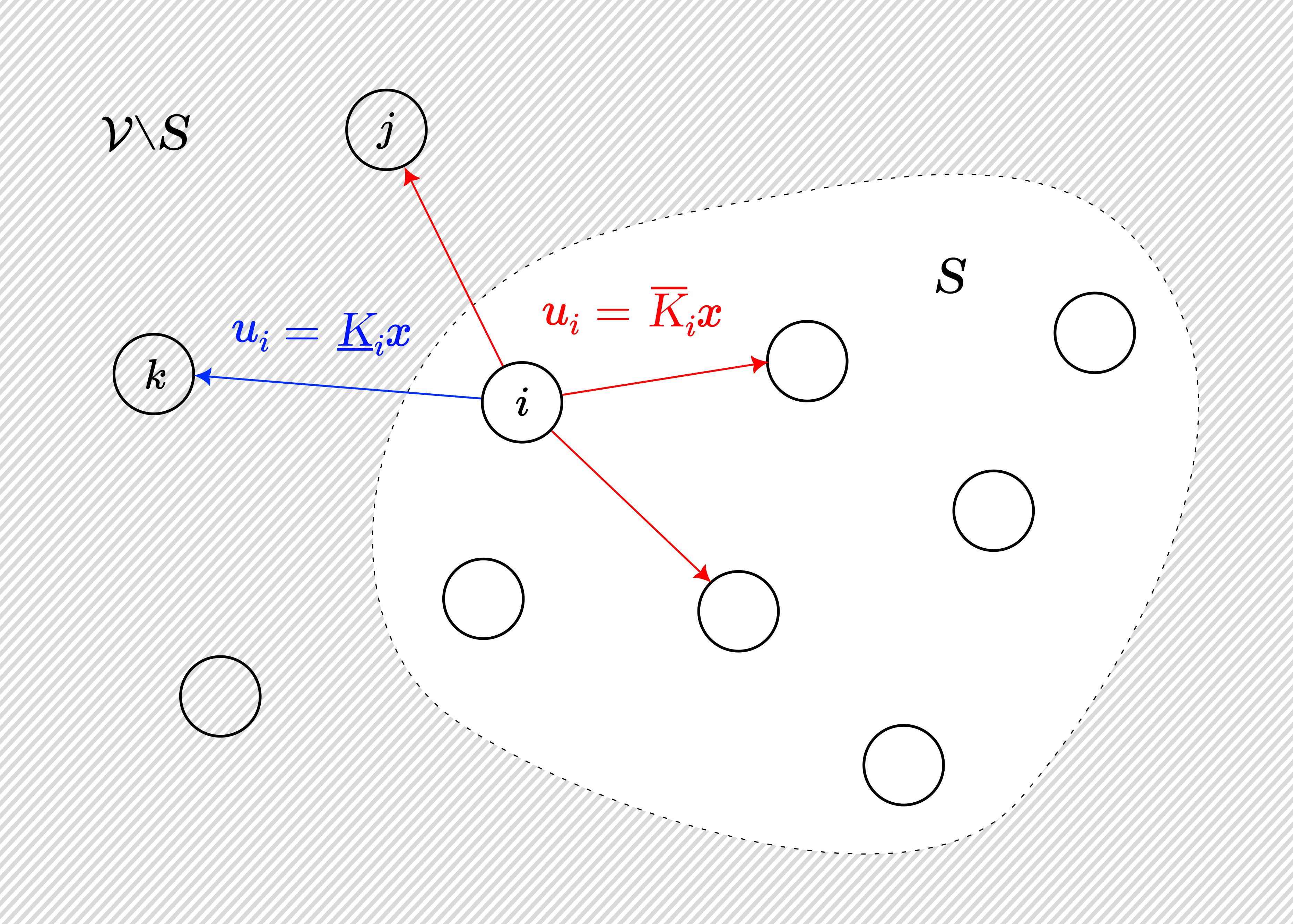}
    \caption{Illustration of the search space $S$ in Algorithm~\ref{alg:1}. Optimal cost functions $\overline{g}_S$ and $\underline{g}_S$ with regard to the control $K_S$ are found using $\overline{h}$ and $\underline{h}$ for states outside~$S$. The space is expanded in each iteration to include states with high uncertainty $\overline{h}_j-\underline{h}_j$ that are most affected under the optimal control laws $\overline{K}_S$ and $\underline{K}_S$ corresponding to $\overline{g}_S$ and $\underline{g}_S$.}
    \label{fig:search}
\end{figure}

\subsection{Stabilized systems}

{In the following setup, we assume knowledge of a stabilizing policy $\hat{K}$ for the full system. This assumption is lifted in the subsequent section.} Let the {subscript} notation {$v_i$} for vectors $v$ denote its $i$th element{, with clarifying parentheses added where subscripts collide (i.e. $(x_0)_i$ indicating the $i$th element of the initial state)}. For the matrices defining the dynamics we define the local dynamics $A_S\in\mathbb{R}^{n\times n}$ as the matrix such that the columns with index $i\in \mathcal{V}\setminus S$ are replaced by unit vectors $e_i$. Further, define two sets of fully actuating controllers on the states in $S$:
\begin{align*}
    \overline{\mathcal{K}}_S &= \{K \in \mathcal{K} : K_i = \hat{K}_i \;\; \forall i\in \mathcal{V}\setminus S\}\\
    \underline{\mathcal{K}}_S &= \{K \in \mathcal{K} : K_i = \mathbf{0} \;\; \forall i\in \mathcal{V}\setminus S\}.
\end{align*} 
These sets represent the use of the stabilizing initial policy $\hat{K}$ outside the search area to guarantee stability in the worst case, while the lower bound uses the optimistic policy ${K_i = \mathbf{0}}$ for states outside $S$.

The resulting dynamics for the search area $S$ are 
\begin{equation*}
    x(t+1) = (A_S + BK)x(t), \;\;\; \overline{\mathcal{K}}_S\;\textnormal{or}\;\underline{\mathcal{K}}_S.
\end{equation*}
This construction preserves any quantity leaving the search area $S$ by transfer to some state $x_i$, $i\in\mathcal{V}\setminus S$ in that state, no longer affecting the states in~$S$. {The final distribution is then used in the prioritization of states when expanding the search area.}

Let $I_S$ denote a matrix with zeros in all elements except diagonal entries 1 at indices corresponding to states in $S$. Define the upper bound local optimum $\overline{g}_S\in\mathbb{R}^{|S|}_+$ as the solution to
\begin{equation}\label{eq:g}
    \overline{g}_S = I_Ss + A_S^\top \overline{g}_S + I_S\min_{K\in{\overline{\mathcal{K}}_S}}\left[K^\top (r + B^\top\overline{g}_S)\right]
\end{equation}
with ${(\overline{g}_S)_i = \overline{h}_i}$ for ${i\in\mathcal{V}\setminus S}$. The lower bound local optimum~$\underline{g}_S$ is defined analogously using ${(\underline{g}_S)_i = \underline{h}_i}$ and ${K\in\overline{\mathcal{K}}_S}$. With these definition we are ready to state {Algorithm~\ref{alg:1}}. The algorithm terminates when the bounds provided by the heuristic function can guarantee that the worst-case cost $\overline{g}^\top x_0$ is at most $\gamma$ times larger than the optimal~$p^\top x_0$.

\begin{algorithm}
\caption{Heuristic search with initial policy}\label{alg:1}
\textbf{Inputs}: $x_0$, $A$, $B$, $E$, $s$, $r$, $\gamma$, $\hat{K}$
\begin{algorithmic}[1]
    \State define $S \leftarrow \{i : (x_0)_i > 0\}$ 
    \State set $\underline{g}_S\leftarrow s$, $\overline{g}_S\leftarrow (I - (A + B\hat{K}))^{-\top}(s+\hat{K}^\top r)$
    \While{$\overline{g}_S^\top x_0/\underline{g}_S^\top x_0 > \gamma$}
    \State find $\overline{g}_S$ and $\underline{g}_S$ that solve \eqref{eq:g} \label{step:p}
    \State\label{step:K} extract {$\overline{K}_S = \textnormal{arg}\min_{K\in\overline{\mathcal{K}}_S} K^\top(r+B^\top\overline{g}_S)$} 
    \State\label{step:x} calculate {${\overline{x} = \lim_{k\rightarrow\infty} (A_S+B\overline{K}_S)^kx_0}$} 
    \State set $S \leftarrow S\cup \textnormal{arg}\max\limits_{i\in\mathcal{V}\setminus S}\left[(\overline{h}_i - \underline{h}_i){\overline{x}_i} \right]$\label{step:S}
    \EndWhile
\end{algorithmic}
\end{algorithm}

\begin{proposition}\label{prop:terminate}
    {Consider a problem on the form \eqref{eq:optprob} with finite value for $x_0\in\mathbb{R}^n_+$.} Algorithm~\ref{alg:1} with $\gamma = 1$ will terminate in finite time and return the optimal cost function in all states for which $x_0$ is non-zero.
\end{proposition}
\vspace*{-4mm}
\begin{pf}
    We first prove that the algorithm terminates. Note that step \ref{step:p} can be solved by means of linear programming as a direct result of Theorem \ref{th:main}, and can thus be completed in finite time. With a view towards generating a contradiction, assume that the algorithm does not terminate. By iteration of step~\ref{step:S}, we must then eventually have $S = \mathcal{V}$. This means that step~\ref{step:p} solves~\eqref{eq:g} for the full state, which is equivalent to~\eqref{eq:p}, yielding ${\overline{g}_S = \underline{g}_S = p}$ according to Theorem~\ref{th:main}{, since the full problem has a finite optimal value}. This in turn means that ${\overline{g}_S^\top x_0/\underline{g}_S^\top x_0 = 1}$, and the algorithm will terminate in the next step, contradicting our assumption.

    The above reasoning shows that the algorithm terminates and that the optimal cost is obtained for all states if ${S = \mathcal{V}}$ upon termination. We now examine the case of termination with ${S\neq\mathcal{V}}$. With ${\gamma = 1}$ this implies that
    \begin{align}\label{eq:gineq}
        \overline{g}_S^\top x_0/\underline{g}_S^\top x_0 &\le 1.
    \end{align}
    Initialization of $\overline{g}_S$ and $\underline{g}_S$ in accordance with Proposition \ref{prop:h} means that ${\underline{g}_S\le p \le\overline{g}_S}$ holds initially. This relation is also preserved under iteration of Algorithm~\ref{alg:1}, as the operation in step \ref{step:p} solves \eqref{eq:g}, which has a strictly larger (in the case of $\overline{g}_S$) solution than the optimal cost of the full system, meaning that ${p\le\overline{g}_S}$ for all iterates. A symmetrical result holds for $\underline{g}_S$. Since ${x_0\ge0}$, the only possible solutions fulfilling~\eqref{eq:gineq} are those for which ${\underline{g}_S^\top x_0 = p^\top x_0 = \overline{g}_S^\top x_0}$. Consequently, we can only have ${(\overline{g}_S)_i \neq p_i}$ and ${(\underline{g}_S)_i \neq p_i}$ for elements $i$ such that ${(x_0)_i = 0}$. \hfill $\blacksquare$
\end{pf}

\begin{corollary}
    Algorithm \ref{alg:1} with parameter ${\gamma\ge1}$ terminates in finite time with cost $\overline{g}^\top x_0 < \gamma p^\top x_0$.
\end{corollary}
\vspace*{-4mm}
\begin{pf}
    For $\gamma > 1$, it is strictly easier to fulfill the termination condition of Algorithm \ref{alg:1}. Thus termination follows directly from Proposition \ref{prop:terminate}. The performance guarantee can be seen from the complement of the termination condition, yielding
    \begin{align*}
        \overline{g}_S^\top x_0/\underline{g}_S^\top x_0 &\le \gamma\\
        \iff \overline{g}_S^\top x_0 &\le \gamma p^\top x_0
    \end{align*}
    where the cost is attained by following the worst-case policy associated with $\overline{g}_S$.
    \hfill $\blacksquare$
\end{pf}
\vspace*{-4mm}

Step \ref{step:p} of Algorithm \ref{alg:1} can, for example, be performed by fixed point iteration of \eqref{eq:g}, or by solving a linear program like ($iii$) in Theorem \ref{th:main}. Either way, the complexity of the calculation scales with the size of the local search area $S$, as opposed to the full state space. In step \ref{step:K}, only local information of $r$ and $B$ is required, as all other elements are multiplied by zero due to the definition of the set $\mathcal{K}_S$. Additionally, if step \ref{step:p} is performed by e.g. value iteration or policy iteration, {the policy $\overline{K}_S$ is} naturally obtained, removing the need to perform step \ref{step:K}.

{As a consequence of Proposition \ref{prop:consistency}, we can show that each iteration and expansion of the search space $S$ monotonically improves the cost.} Meanwhile, as noted previously, the conditions for a consistent heuristic are identical to the constraints of the linear program. {This means that, when solving step~\ref{step:p} using linear programming, the previous solution is a feasible point from which we can initialize the program of the next iteration.} 
\vspace*{-2mm}
{\begin{proposition}
    Given a stabilizing policy $\hat{K}$, Algorithm~\ref{alg:1} converges monotonically.
\end{proposition}}
\vspace*{-4mm}
{\begin{pf}
    Let $S_1$ and $S_2$ be the sets corresponding to the search area in subsequent iterations of Algorithm~\ref{alg:1}. Consequently, $S_1\subset S_2$. When evaluating~\eqref{eq:g} for $S_2$, we see that the contribution to $\overline{g}_{S_2}$ from states outside $S_2$ is identical to the case of~\eqref{eq:g} for $S_1$. The only difference, then, is the influence of the state, $i$ say, added to the search area between the two iterations. This influence is chosen as $\min\{r_i+B_i^\top\overline{g}_{S_2}\}E_i$, rather than being fixed to the previous suboptimal policy. We now have that ${\overline{g}_{S_1}\le\overline{h}}$ implies ${\overline{g}_{S_2}\le\overline{g}_{S_1}}$. The first step of the induction follows immediately from consistency of $\overline{h}$, as given in Proposition~\ref{prop:h}. \hfill $\blacksquare$
\end{pf}}
\vspace*{-5mm}

We note that the optimization over local controllers in step~\ref{step:K} can be simplified in the event that the expression 
\begin{equation*}
    \textnormal{arg}\min_{K\in\mathcal{K}_i} K^\top(r_i+B_i^\top g)
\end{equation*}
for some index $i$ yields the same result for all $g$ in the interval $\underline{g}_S \le g \le \overline{g}_S$. In this case, one action is optimal for all possible values of $p$ within the current bounds. This action can thus be fixed in all future iterations, reducing the computational complexity. This idea of "locking in" actions that can be certified as optimal within the given bounds is leveraged in several heuristic search frameworks for SSP, see e.g. the Hierarchical Dirichlet Process (HDP) algorithm in \citep{bonet03faster}.

Step \ref{step:x} is significantly reduced in size for systems with sparse dynamics. The state $\overline{x}$ is only calculated for states that are directly impacted by those in the search area $S$. The priority of states outside $S$ in step \ref{step:S} is not uniquely optimal, see \citep{mcmahan05brtdp} for a discussion of different methods. The current setting differs in that dynamics in continuous state are simulated to determine which states have the most potential to improve the total cost, both as a consequence of a large uncertainty $\overline{h_i}-\underline{h}_i$ and as a result of a high degree of utilization $\overline{x}_i$ in the worst case. Changing or weighting this decision offers the chance to tailor the algorithm for, e.g., more optimistic or pessimistic intermediate solutions.
\vspace*{-2mm}
{\begin{remark}
    In the case $E=I$, the heuristics contain all information required to find the optimal policy for a subset $S$ of nodes. In other words, we can find the optimal parameter of~\eqref{eq:p} by solving the smaller subproblem~\eqref{eq:g} with local policy $K\in\overline{\mathcal{K}}_S$ when the optimal solution is known for states outside $S$. This is no longer the case when $E\neq I$, as the sum in the expression~\eqref{eq:p} then indicates that the optimal policy in states outside $S$ influence the optimal value. As a consequence, we expect deteriorated performance of the heuristic search for $E$ with large off-diagonal elements. This serves to illustrate the difference between the class of SSP and the full class~\eqref{eq:optprob}, and explains the remarkable success of heuristic methods for SSP in particular.
\end{remark}}

\subsection{{Unstable systems}}

{In situations where no initial stabilizing policy is known, heuristic bounds can still help guide the search for an optimal policy. However, the method presented above must be modified. Consider the following implicit equation for the upper bound solution on the search area $S$:}
\begin{equation}\label{eq:newg}
    {\overline{g}_S = I_Ss + A_S^\top \overline{g}_S + I_S\min_{K\in{\underline{\mathcal{K}}_S}}\left[K^\top (r + B^\top\overline{g}_S)\right]}
\end{equation}
{with $\overline{g}_i = \overline{h}_i$ for $i\in\mathcal{V}\setminus S$. The right-hand expression is identical to that in~\eqref{eq:g} apart from the use of the zero policy $K_i = \mathbf{0}$ outside the search area. This represents the lack of information about the behavior of the system outside $S$, removing the influence of the unknown inputs $u_i$ for $i\in\mathcal{V}\setminus S$ on the nodes in $S$. As in the previous section, an equivalent expression for $\underline{g}_S$ is obtained by instead using the lower bound heuristic $(\underline{g}_S)_i = \underline{h}_i$ for states $i$ outside $S$. The procedure for heuristic search without prior knowledge of a stabilizing policy is presented in Algorithm~\ref{alg:unstable}.} {Just as above,~\eqref{eq:newg} coincides with~\eqref{eq:p} for $S = \mathcal{V}$, guaranteeing convergence of Algorithm~\ref{alg:unstable} within the prescribed performance bound~$\gamma$. This follows directly from the same reasoning as used in the proof of Proposition~\ref{prop:terminate} and Corollary~1.}

\begin{algorithm}
\caption{{Heuristic search without initial policy}}\label{alg:unstable}
\textbf{Inputs}: $x_0$, $A$, $B$, $E$, $s$, $r$, {$\overline{h}$, $\underline{h}$}, $\gamma$
\begin{algorithmic}[1]
    \State define $S \leftarrow \{i : (x_0)_i > 0\}$ 
    \State set ${\underline{g}\leftarrow \underline{h}}$, ${\overline{g}\leftarrow \overline{h}}$
    \While{$\overline{g}_S^\top x_0/\underline{g}_S^\top x_0 > \gamma$}
    \State find $\overline{g}_S$ and $\underline{g}_S$ that solve {\eqref{eq:newg}} \label{step:p}
    \State\label{step:K} extract {$\overline{K}_S = \textnormal{arg}\min_{K\in\underline{\mathcal{K}}_S} K^\top(r+B^\top\overline{g}_S)$} 
    \State\label{step:x} calculate {${\overline{x} = \lim_{k\rightarrow\infty} (A_S+B\overline{K}_S)^kx_0}$} 
    \State set $S \leftarrow S\cup \textnormal{arg}\max\limits_{i\in\mathcal{V}\setminus S}\left[(\overline{h}_i - \underline{h}_i){\overline{x}_i} \right]$\label{step:S}
    \EndWhile
\end{algorithmic}
\end{algorithm}

{Special care needs to be taken when evaluating the expansion of $S$ in step \ref{step:S} of Algorithm~\ref{alg:unstable}. Since $\overline{h}$ is not guaranteed to be finite for all states, we let the product $(\overline{h}_i - \underline{h}_i)\overline{x}_i$ take the value 0 for states where $\overline{x}_i=0$. This means that the algorithm prioritizes exploring unstable states that cannot be avoided in the worst case $\overline{x}$. Unstable states that can be avoided (giving $\overline{x}_i = 0$) are deprioritized, as the heuristics give no information to guide the search in these states.} {Note that this implies that a stabilizing policy for the states in $S$ has been found when Algorithm~\ref{alg:unstable} terminates, regardless of $\gamma$.} 
\vspace*{-2mm}
{\begin{proposition}
    Given admissible heuristic bounds $\overline{h}$ and $\underline{h}$ such that $\underline{h}$ is consistent and $\overline{h}$ is superconsistent, Algorithm~\ref{alg:unstable} converges monotonically.
\end{proposition}}
\vspace*{-4mm}
{\begin{pf}
    Let $S_1$ and $S_2$ be the sets corresponding to the search area in subsequent iterations of Algorithm~\ref{alg:unstable}, and let $i$ be the index of the state added between the two iterations, so that $S_2=S_1\cup \{i\}$. We wish to show that $\overline{g}_{S_1}\ge\overline{g}_{S_2}$. Since $\overline{h}$ is superconsistent, the element corresponding to the added state fulfills
    \begin{align*}
        (\overline{g}_{S_1})_i &= \overline{h}_i\\
        &\ge s_i + A_i^\top \overline{h} + \min\{r_i+B_i^\top\overline{h},0\}E_{ii}\\
        &\ge s_i + A_i^\top \overline{h} + \sum_{j=1}^n\min\{r_j+B_j^\top\overline{g}_{S_1},0\}E_{ji}\\
        &= (\overline{g}_{S_2})_i.
    \end{align*}
     Here, the second inequality follows admissibility of $\overline{h}$ and the fact that all terms of the right-hand side sum are nonpositive.
    \hfill $\blacksquare$
\end{pf}}
\vspace*{-10mm}
{\begin{remark}
    Note that, in order to attain the specified performance bound and guarantee stability, we implicitly assume that performance within the heuristic bounds $\overline{h}$ and $\underline{h}$ can be realized. This does not mean, however, that we necessarily know the policy in states outside $S$. The heuristics may be obtained from measurements of performance under some nominal operation of the system.
\end{remark}}

\section{{Distributed value iteration}}

{Algorithms~\ref{alg:1} and~\ref{alg:unstable} leverage the constructed heuristic bounds to find close-to-optimal solutions for specific initializations $x_0$. To achieve this, centralized information about the dynamics is required for the growing search region. Below, we present a second algorithm that instead uses the local bound on performance in each state implied by the heuristics to distribute the value iteration algorithm. In contrast to standard Gauss-Seidel iteration used to distribute value iteration of \eqref{eq:p} (as used in \citep{ohlin23optimal}), the use of upper and lower bounds allows for a distributed stopping condition based on the desired distance to optimality, with global performance guarantees.}

{Consider $n$ agents, each representing one element $x_i$ of the state vector. In the graphical interpretation of Figure~\ref{fig:search}, each agent corresponds to a node, locally determining the gain $K_i$ for the inputs $u_i$ and corresponding constraint. Each matrix $B_i$ as defined in the problem setup (Section~2) above governs the effect of agent $i$'s local actions. Define the local incidence matrix ${W^{(i)} = B_iB_i^\top}$ for each agent, and the following sets of neighbors: }
\begin{align*}
    &{\mathcal{N}^{A}_i := \{j: A_{ji}\neq 0, j\neq i\}}\\
    &{\mathcal{N}^E_i = \{j: E_{ji}\neq 0, j\neq i\}}\\
    &{\mathcal{N}_i^{B} = \{j: W^{(i)}_{ij}\neq 0, j\neq i\}.}
\end{align*}

{Algorithm \ref{alg:2} successively improves the local bounds on the cost function for each agent via value iteration. Agent~$i$ terminates when the local suboptimality gap is reduced below a chosen parameter $\gamma$. The global performance of the algorithm is shown in the following proposition. It proves that termination by one or more agents does not prevent the remainder from reaching the prescribed level of optimality. In other words, the suboptimality gap induced by agents terminating does not grow beyond the factor $\gamma$ as it propagates through the network, regardless of the number of agents. This performance is achieved by the resulting pessimistic policy $\overline{K}$, associated with the obtained upper bound.}

\begin{algorithm}
{
\caption{{Distributed value iteration}}\label{alg:2}
\textbf{Inputs}: $A$, $B$, $E$, $s$, $r$, $\overline{h}$, $\underline{h}$, $\gamma$
\begin{algorithmic}[1]
    \State define $\mathcal{I} = \{1,...,n\}$
    \ForAll{$i\in\mathcal{I}$}
    \State set $\overline{g}_i\leftarrow \overline{h}_i$, $\underline{g}_i\leftarrow \underline{h}_i$, $\overline{q}_i\leftarrow 0$, $\underline{q}_i\leftarrow 0$
    \EndFor
    \While{$\mathcal{I}\neq \emptyset$}
    \State sample agent $i\in \mathcal{I}$
    \State receive $\overline{g}_j$, $\underline{g}_j$ from each neighbor $j \in \mathcal{N}_i^A \cup \mathcal{N}_i^B$
    \State receive $\overline{q}_j$, $\underline{q}_j$ from each neighbor $j \in \mathcal{N}_i^E$
    \State $\overline{q}_i \gets \min\{r_i + B_i^\top \overline{g}, 0\}$, $\underline{q}_i \gets \min\{r_i + B_i^\top \underline{g}, 0\}$\label{step:q}
    \State $\overline{g}_i \gets s_i + A_i^\top \overline{g} + \overline{q}_iE_{ii} + \sum\limits_{j \in \mathcal{N}_i^E} \overline{q}_jE_{ji}$
    \State $\underline{g}_i \gets s_i + A_i^\top \underline{g} + \underline{q}_iE_{ii} + \sum\limits_{j \in \mathcal{N}_i^E} \underline{q}_jE_{ji}$
    \If{$\overline{g}_i/\underline{g}_i <\gamma$}
    \State set $\mathcal{I} \leftarrow \mathcal{I}\setminus i$
    \EndIf
    \EndWhile
\end{algorithmic}}
\end{algorithm}

\begin{proposition}
    {Consider a problem on the form \eqref{eq:optprob} with finite value for $x_0\in\mathbb{R}^n_+$. Algorithm \ref{alg:2} converges to a solution $\overline{g}$ such that $\overline{g} \le \gamma p$ for any given $\gamma \ge 1$.}
\end{proposition}
\vspace*{-4mm}
\begin{pf}
    {Algorithm~\ref{alg:2} performs value iteration for $\overline{g}$ and $\underline{g}$ using the Gauss-Seidel method, which converges asymptotically to the fixed point of \eqref{eq:p}. When an agent terminates, its current (suboptimal) bounds are fixed. This means that the iterations of upper and lower bounds for the remaining system have two separate fixed points $\overline{g}^*$ and $\underline{g}^*$. To prove the stated performance bound, we will show that termination of one or more agents under the given condition cannot perturb the distance between the fixed points of the upper and lower bounds for the remaining system beyond a factor $\gamma$.}

    {
    Consider the scenario where a subset $\mathcal{I}$ of agents are still iterating according to Algorithm \ref{alg:2}, while the rest have terminated, partitioning the upper bound vector accordingly as $\overline{g}_{_\mathcal{I}}$ and $\overline{g}_{_{-\mathcal{I}}}$. Rewriting \eqref{eq:Kform} gives the following expression for the fixed point:}
    \begin{equation}\label{eq:M}{
        \begin{bmatrix}
            \overline{g}^*_{_\mathcal{I}} \\ \overline{g}^*_{_{-\mathcal{I}}}
        \end{bmatrix} = \underbrace{\begin{bmatrix}
            \overline{M}_{_{\mathcal{I},\mathcal{I}}} & \overline{M}_{_{\mathcal{I},-\mathcal{I}}} \\ \overline{M}_{_{-\mathcal{I},\mathcal{I}}} & \overline{M}_{_{-\mathcal{I},-\mathcal{I}}}
        \end{bmatrix}}_{= (A + B\overline{K})^\top} \begin{bmatrix}
            \overline{g}^*_{_\mathcal{I}} \\ \overline{g}^*_{_{-\mathcal{I}}}
        \end{bmatrix} + \underbrace{\begin{bmatrix}
            \overline{k}_{_\mathcal{I}} \\ \overline{k}_{_{-\mathcal{I}}}
        \end{bmatrix}}_{= s + \overline{K}^{\top} r}.
        }
    \end{equation}
    {Solving for $\overline{g}^*_{_\mathcal{I}}$ we get}
    \begin{equation}\label{eq:fixed}{
        \overline{g}^*_{_\mathcal{I}} = (I-\overline{M}_{_{\mathcal{I},\mathcal{I}}})^{-1} \left(\overline{M}_{_{\mathcal{I},-\mathcal{I}}}\overline{g}^*_{_{-\mathcal{I}}} + \overline{k}_{_{\mathcal{I}}}\right).}
    \end{equation}
    {An equivalent expression holds for $\underline{g}^*_{_\mathcal{I}}$, with dynamics~$\underline{M} = (A+B\underline{K})$. 
    In order to prove that the difference between the perturbed fixed points is still within the specified bound, we need to show nonpositivity of the quantity}
    \begin{multline}\label{eq:gdiff}
        {\overline{g}^*_{_\mathcal{I}}-\gamma\underline{g}^*_{_\mathcal{I}} = \overline{g}^*_{_{\mathcal{I}}}- \gamma (I-\underline{M}_{_{\mathcal{I},\mathcal{I}}})^{-1} \left(\underline{M}_{_{\mathcal{I},-\mathcal{I}}}\underline{g}^*_{_{-\mathcal{I}}} + \underline{k}_{_{\mathcal{I}}}\right).}
    \end{multline}
    {For agents not in $\mathcal{I}$, the stopping condition implies that ${\overline{g}^*_{_{-\mathcal{I}}} \le \gamma \underline{g}^*_{_{-\mathcal{I}}}}$. Further, the matrix $\underline{M}_{_{\mathcal{I},\mathcal{I}}}$ is a principal submatrix of $(A+B\underline{K})^\top$ and thus Schur stable. This together with $\underline{M}_{_{\mathcal{I},\mathcal{I}}}\ge 0$ means that $I-\underline{M}_{_{\mathcal{I},\mathcal{I}}}$ is an M-matrix, which in turn implies that $(I-\underline{M}_{_{\mathcal{I},\mathcal{I}}})^{-1}\ge0$. We therefore obtain the following element-wise upper bound on \eqref{eq:gdiff}:} 
    \begin{multline}\label{eq:upper}
        {\overline{g}^*_{_\mathcal{I}}-\gamma\underline{g}^*_{_\mathcal{I}} \le \overline{g}^*_{_{\mathcal{I}}}- (I-\underline{M}_{_{\mathcal{I},\mathcal{I}}})^{-1} \left(\underline{M}_{_{\mathcal{I},-\mathcal{I}}}\overline{g}^*_{_{-\mathcal{I}}} + \underline{k}_{_{\mathcal{I}}}\right)}\\ {-(\gamma-1)(I-\underline{M}_{_{\mathcal{I},\mathcal{I}}})^{-1} \underline{k}_{_{\mathcal{I}}}.}
    \end{multline}
    {Here, the terms have been rearranged to illustrate the fact that the second term corresponds to the expression~\eqref{eq:fixed} for the upper bound if the optimal policy is replaced by $\underline{K}$. This policy cannot outperform~
    $\overline{K}$, and Proposition~\ref{prop:consistency} shows that the optimal upper bound is \textit{element-wise} less than or equal to those under other policies. Thus, the sum of the first two terms in \eqref{eq:upper} is nonpositive in every element. Nonpositivity of the right-hand side in~\eqref{eq:upper} now follows from $\underline{k}_{_{\mathcal{I}}} > 0$, completing the proof.}\hfill $\blacksquare$
\end{pf}

\begin{figure*}[t]
    \centering
    \includegraphics[width=\textwidth]{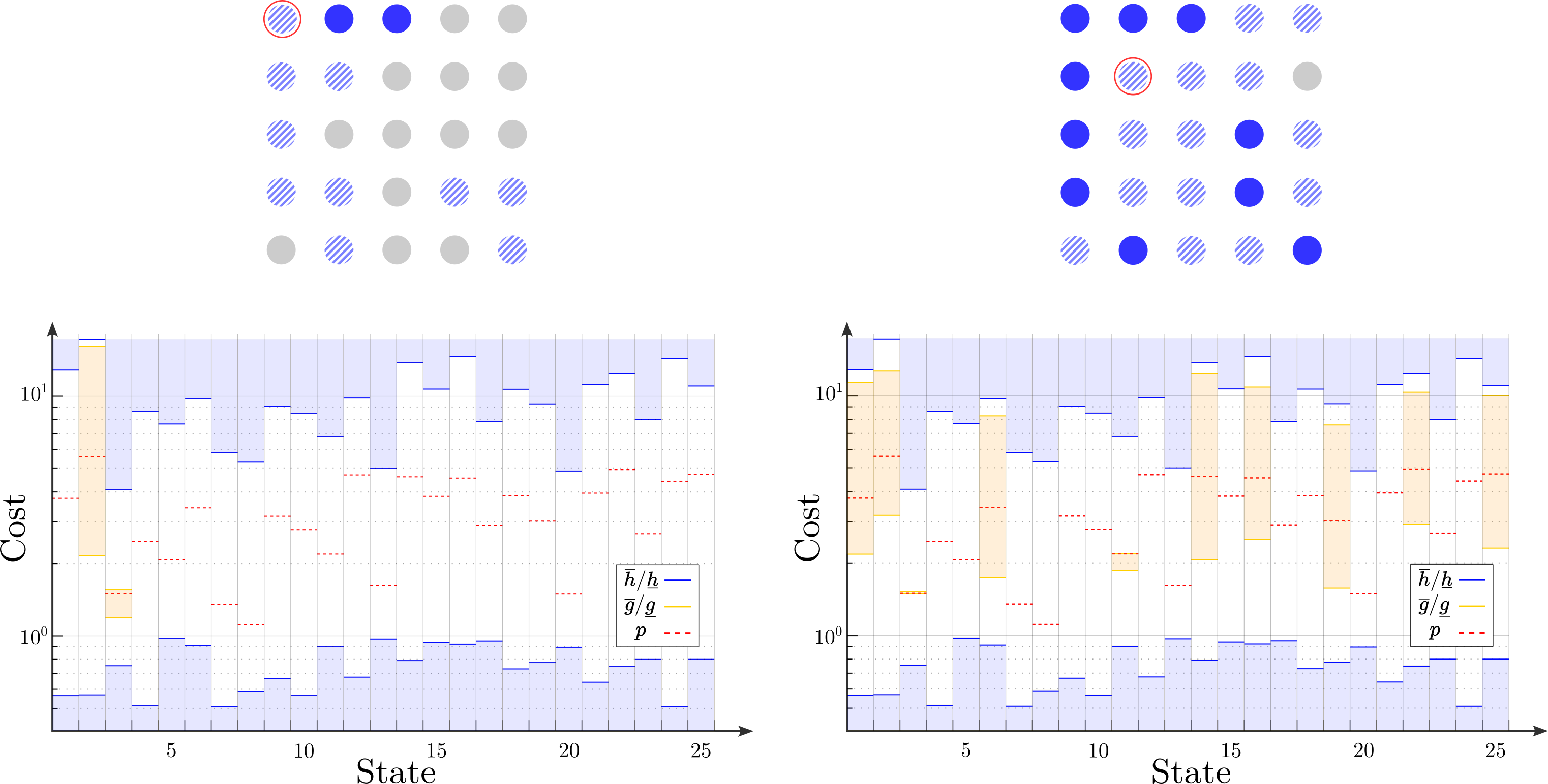}
    \caption{Illustration of the search space after the first (left) and ninth (right) iteration of Algorithm~\ref{alg:1} for the system in Example~\ref{ex:big}. Above, states in dark blue are part of the search space $S$ while shaded states are in the active boundary of states that can be directly affected by the dynamics inside $S$. The red circle indicates the state selected for inclusion in the next iteration. The accompanying plots show the optimal cost function~$p$ (dashed red) for each state, as well as the current estimates $\overline{g}$ and $\underline{g}$ (yellow) for states in~$S$ and the heuristic bounds $\overline{h}$ and $\underline{h}$ (blue). {It is clear from the schematic view of the state space that the algorithm systematically explores states most impacted from the initial configuration. Note in particular that the performance of $x_3$ is quickly controlled to within the performance bound, while $x_2$ remains uncertain, necessitating the exploration of a large part of the remaining state space.}}
    \label{fig:searchplot}
\end{figure*}

\begin{figure}[t]
    \centering
    \includegraphics[width=.9\linewidth]{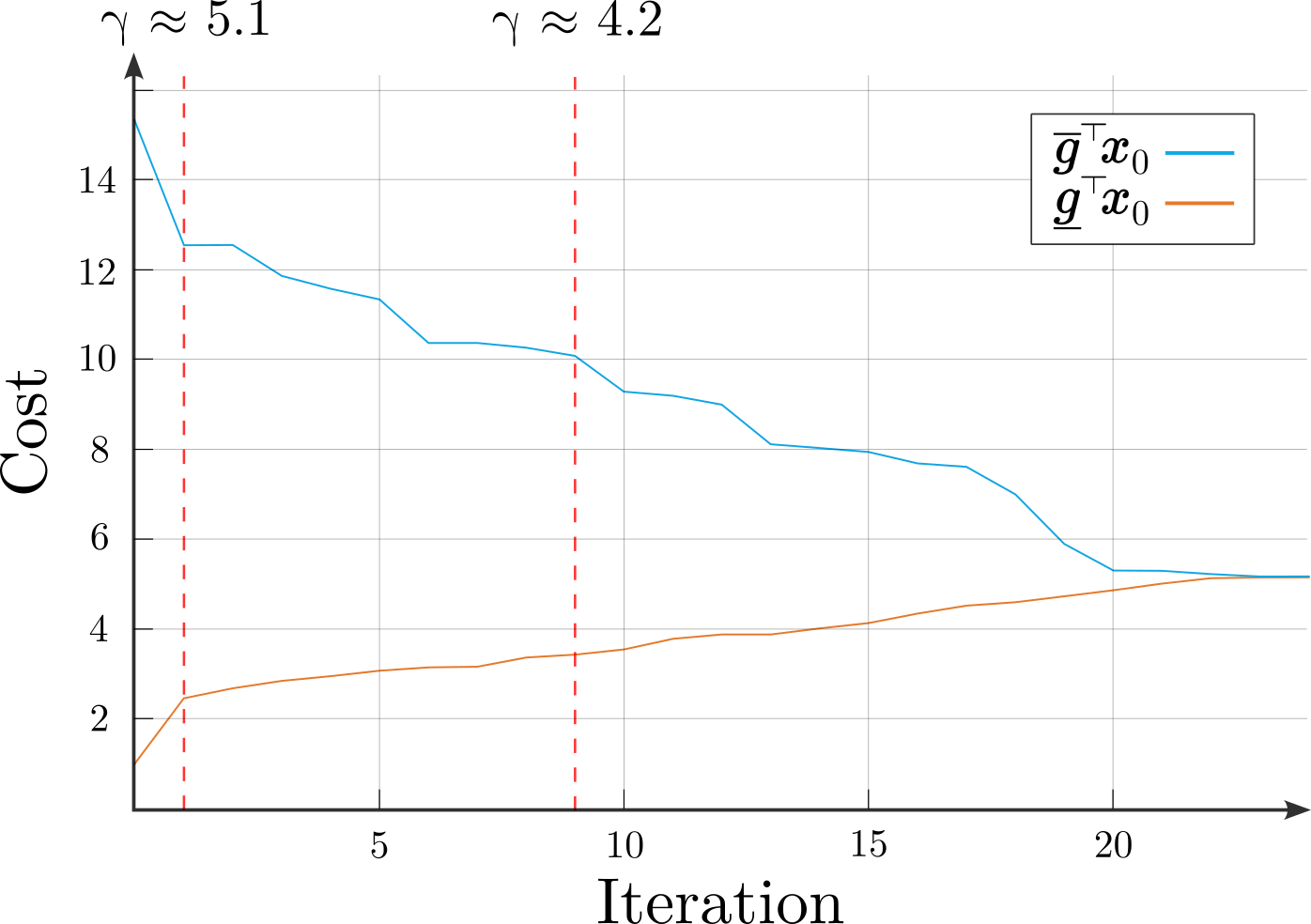}
    \caption{Upper and lower bound of the total cost in Example~\ref{ex:big}. The red lines highlight the iterations for which the state is displayed in more detail in Figure~\ref{fig:searchplot}. The value~$\gamma$ specifies the stopping condition in Algorithm~\ref{alg:1} and guarantees a performance of at least $\overline{g}^\top x_0 \le \gamma p^\top x_0$ when the algorithm terminates. {Since the initial heuristics in Example~\ref{ex:big} are consistent, the convergence is monotone.}}
    \label{fig:cost}
\end{figure}

\section{Numerical example}

Next, we give a simulated example to illustrate the properties of Algorithm~\ref{alg:1}.

\begin{example}\label{ex:big}

\textnormal{To give a plausible example of a process that satisfies the necessary assumptions, we consider a large chemical plant. The state $x\in\mathbb{R}^{25}$ represents amounts of $n = 25$ different chemical waste compounds that are to be disposed at as low cost (corresponding to energy consumption, pollution or resources consumed) as possible. For each compound a variety of reactions can be induced, represented by the dynamics $Ax+Bu$, which give the resultant compounds as the updated state. We generate a random matrix ${A\ge0}$ with some sparsity and construct ${B\in\mathbb{R}^{25\times75}}$ such that ${E=A}$ guarantees positivity of the dynamics, see \citep{ohlin23optimal}.}

\textnormal{Let $m_i = 2$ for all $i$ giving a total of three reactions, one represented by the autonomous dynamics and the other two by inputs $u_{ij}$. The cost of each reaction is modeled by vectors $s$ and $r$. Further, the first action for each state $u_{i1}$ is modeled as a direct disposal of all compound of the corresponding state at a high cost. This could represent e.g. burning or releasing the dangerous chemical, resulting in large amounts of pollution. In the input matrix $B$ this is modeled as the column multiplying $u_{i1}$ having its only non-zero element $-1$ on its $i$th row. 
The initial stabilizing policy $\hat{K}$ is constructed by taking the action $u_{i1}$ for all $i$, giving a suboptimally high cost but guaranteeing that the state is stable. As a consequence, $A+B\hat{K} = \mathbf{0}$. The expression for the upper bound heuristic in Proposition \ref{prop:h} then becomes}
\begin{equation*}
    \overline{h} = s + \hat{K}^\top r.
\end{equation*}
\textnormal{Note that, due to the above construction of $B$, this is equivalent to $\overline{h} = s + A^\top\hat{r}$ where $\hat{r}\in\mathbb{R}^{25}$ is the subvector of entries of~$r$ corresponding to the initial policy. As the lower bound policy is given by $\underline{h} = s$, the initial uncertainty in the cost function is quantified by~$A^\top\hat{r}$. }

\textnormal{We run Algorithm~\ref{alg:1} starting from the initial state ${x_2 = 0.7}$, ${x_3 = 0.8}$ and ${x_i = 0}$ for all other states $i$. Figure \ref{fig:cost} displays the upper and lower bound for the cost incurred starting from $x_0$, as a function of the size of the search space $S$. Due to the connected nature of the dynamics and the poor performance of the initial policy, fully optimal control necessitates exploring the full state space. In Figure \ref{fig:searchplot}, two snapshots (indicated in Figure \ref{fig:cost}) of the algorithm are shown. Each illustrates the current search space, as well as the values of the intermediate bounds $\overline{g}$ and $\underline{g}$ compared to the initial heuristic bounds and the optimal value $p$.}

\end{example}

\section{Conclusion}

The presented work builds on recent results for positive linear systems with explicit solutions. This allows us to demonstrate an equivalence between a {subset} of such problems and the class of stochastic shortest path problems that are the focus of heuristic search methods. An example of such a method applied to the optimal control setting is given, which gives results for the full class of optimal control problems, broadening the scope of applicability beyond SSP. This connection can serve as a bridge to gain a wider field of application for existing results in both areas. {An algorithm for distributed optimal control is presented, where the derived heuristic bounds are the key tool to enable global bounds on performance.}

We argue for the possibility of designing heuristic search algorithms, with several initial examples provided here, inspired by existing methods for shortest path problems to achieve scalable distributed control of large-scale networked systems. {Promising areas of further study include the dual to the presented problem, which suggests that favorable results could be derived for state observers of positive systems.}

\begin{ack}

We would like to thank Olle Kjellqvist, Jonas Hansson {and Mark Jeeninga} at Lund university {and César Uribe at Rice Academy} for insightful discussion related to the presented material. {We are also grateful to the anonymous reviewers for their constructive comments and helpful advice.}

\end{ack}

\bibliographystyle{agsm}
\bibliography{references}           

\begin{wrapfigure}{l}{20mm} 
    \includegraphics[width=1in,height=1.25in,clip,keepaspectratio]{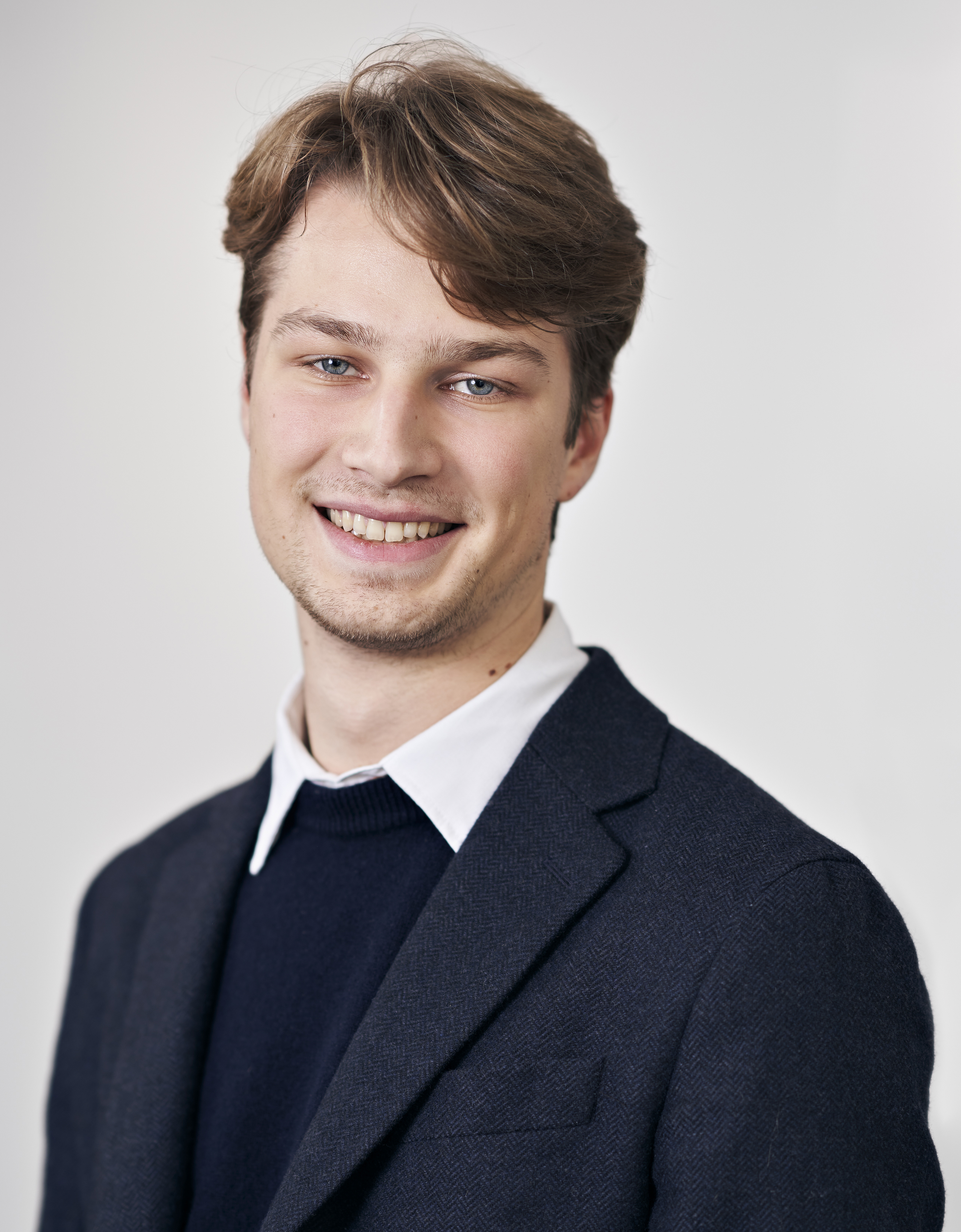}
\end{wrapfigure}\par

\footnotesize{\textbf{David Ohlin} is a Ph.D. student at the Department of Automatic Control at Lund University, Sweden, where he received his M.Sc. in Engineering Physics in 2021. He was a visiting researcher at the Department of Electrical Engineering and Computer Sciences at UC Berkeley in 2024. His research focuses on the modeling of dynamical systems on graphs, with applications in scalable control and opinion dynamics.\par}

\begin{wrapfigure}{l}{20mm} 
    \includegraphics[width=1in,height=1.25in,clip,keepaspectratio]{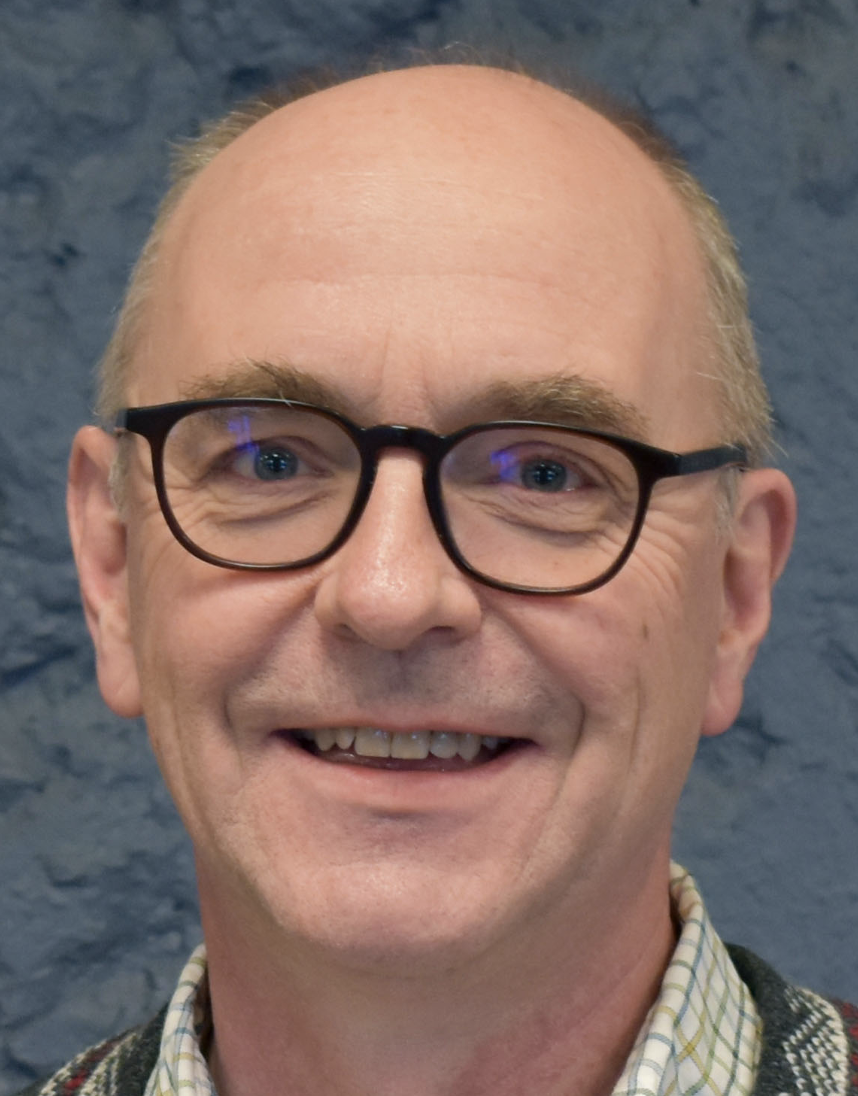}
\end{wrapfigure}\par
\footnotesize{\textbf{Anders Rantzer} was appointed professor of Automatic Control at Lund University, Sweden, after a PhD at KTH Stockholm in 1991 and a postdoc 1992/93 at IMA, University of Minnesota. He was visiting associate faculty member at Caltech the academic year of 2004/05, Taylor Family Distinguished Visiting Professor at University of Minnesota 2015/16 and for 2026-2028 he has been appointed Honorary Chair Professor of National Sun Yat-sen University, Taiwan. Rantzer is elected member of the Royal Swedish Academy of Engineering Sciences and a Fellow of IEEE and IFAC. He has been chairman of the Swedish Scientific Council for Natural and Engineering Sciences as well as the Royal Physiographic Society in Lund. His research interests are in modeling, analysis and synthesis of control systems, with particular attention to scalability and adaptation.\par}

\begin{wrapfigure}{l}{20mm} 
    \includegraphics[width=1in,height=1.25in,clip,keepaspectratio]{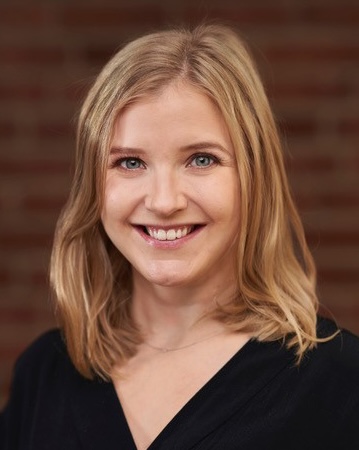}
\end{wrapfigure}

\textbf{Emma Tegling} is an Associate Professor with the Department of Automatic Control at Lund University, Sweden and a fellow of the Wallenberg AI, Autonomous Systems and Software Program (WASP). She received her Ph.D.in Electrical Engineering from KTH Royal Institute of Technology, Stockholm, Sweden in 2019, and her M.Sc. and B.Sc. degrees, both in Engineering Physics, from the same institute in 2013 and 2011, respectively. From 2019 to 2020 she was a Postdoctoral Research Fellow with the Institute of Data, Systems, and Society (IDSS) at the Massachusetts Institute of Technology (MIT), Cambridge, USA. She is a member of the Board of Governors of the IEEE Control Systems Society. Her research interests revolve around large-scale network systems, their structure and control, fundamental limitations, and applications including power networks and social networks. 

\end{document}